\documentclass[reqno]{amsart}
\usepackage{times}
\usepackage{amsmath,amssymb}
\usepackage{graphicx}
\usepackage{color}
\usepackage[latin1]{inputenc}
\usepackage[all]{xy}
\usepackage[dvips]{hyperref}

\newcommand {\NN}{{\mathbb N}}
\newcommand {\PP}{{\mathbb P}}

\newcommand {\RR}{{\mathbb R}}
\newcommand {\ZZ}{{\mathbb Z}}

\DeclareMathOperator {\codim}{codim}

\DeclareMathOperator {\Div}{Div}
\DeclareMathOperator {\Pic}{Pic}

\DeclareMathOperator {\divisor}{div}
\DeclareMathOperator {\id}{id}
\DeclareMathOperator {\aff}{aff}

\DeclareSymbolFont {mysymbols}{OMS}{cmsy}{m}{n}
\DeclareMathSymbol {\calI}{\mathalpha}{mysymbols}{`I}
\DeclareMathSymbol {\calK}{\mathalpha}{mysymbols}{`K}
\DeclareMathSymbol {\calO}{\mathalpha}{mysymbols}{`O}
\DeclareMathSymbol {\calT}{\mathalpha}{mysymbols}{`T}

\DeclareSymbolFont {mylargesymbols}{OMX}{cmex}{m}{n}
\DeclareMathSymbol {\dunion}{\mathop}{mylargesymbols}{"60}

\newcommand {\df}[1]{\textsl {#1}}

\newcommand {\preprint}[2]{preprint \discretionary {#1/}{#2}{#1/#2}}

\newcommand {\twolines}[2]{\stackrel{\scriptstyle {#1}}{#2}}
\newcommand {\subfan}[0]{\!\unlhd}

\renewenvironment {enumerate}%
  {\rule{1mm}{0mm}\begin {oldenumerate}%
    \parskip1ex plus0.5ex \itemsep 0mm \parindent 0mm}%
  {\end {oldenumerate}}

\renewenvironment {itemize}%
  {\rule{1mm}{0mm}\begin {olditemize}%
    \parskip1ex plus0.5ex \itemsep 0mm \parindent 0mm}%
  {\end {olditemize}}

\theoremstyle {plain}
\newtheorem {theorem}{Theorem}[section]
\newtheorem {proposition}[theorem]{Proposition}
\newtheorem {lemma}[theorem]{Lemma}
\newtheorem {corollary}[theorem]{Corollary}
\theoremstyle {definition}
\newtheorem {definition}[theorem]{Definition}
\theoremstyle {remark}
\newtheorem {remark}[theorem]{Remark}
\newtheorem {remdef}[theorem]{Remark and Definition}
\newtheorem {example}[theorem]{Example}
\newtheorem {construction}[theorem]{Construction}

\begin{document}

  \title[First Steps in Tropical Intersection Theory]{First Steps in Tropical Intersection Theory}
\author {Lars Allermann and Johannes Rau}
\address {Lars Allermann, Fachbereich Mathematik, TU Kaiserslautern, Postfach
  3049, 67653 Kaiserslautern, Germany}
\email {allerman@mathematik.uni-kl.de}
\address {Johannes Rau, Fachbereich Mathematik, TU Kaiserslautern, Postfach
  3049, 67653 Kaiserslautern, Germany}
\email {jrau@mathematik.uni-kl.de}

\begin {abstract}
We establish first parts of a tropical intersection theory.
Namely, we define cycles, Cartier divisors and intersection
products between these two ({\itshape without} passing to rational
equivalence) and discuss push-forward and pull-back. We do this
first for fans in $\RR^n$ and then for ``abstract'' cycles that
are fans locally. With regard to applications in enumerative
geometry, we finally have a look at rational equivalence and
intersection products of cycles and cycle classes in $\RR^n$.
\end {abstract}

\maketitle

  \section{Introduction} \label{sec-introduction}

Tropical geometry is a recent development in the field of
algebraic geometry that tries to transform algebro-geometric
problems into easier, purely combinatorial ones. In the last few
years various authors were able to answer questions of enumerative
algebraic geometry using these techniques. In order to determine
the number of (classical) curves meeting given conditions in some
ambient space they constructed moduli spaces of tropical curves
and had to intersect the corresponding tropical conditions in
these moduli spaces. Since there is no tropical intersection
theory yet the computation of the arising intersection
multiplicities and the proof of the independence of the choice of
the conditions had to be repeated for every single problem without
the tools of an elaborated intersection theory (see for example
\cite{GM}, \cite{KM}).

A first draft of a general tropical intersection theory without
proofs has been presented by Mikhalkin in \cite{M}. The concepts
introduced there --- if set up rigorously --- would help to unify
and solve the above mentioned problems and would provide utilities
for further applications. Thus in this paper we develop in detail
the basics of a general tropical intersection theory based on
Mikhalkin's ideas.

This paper consists of three parts: In the first part (sections
\ref{sec-affinecycles} - \ref{sec-pushpull}) we firstly introduce
affine tropical cycles as balanced weighted fans modulo
refinements and affine tropical varieties as affine cycles with
non-negative weights. One would like to define the intersection
of two such objects but in general neither is the set-theoretic
intersection of two cycles again a cycle nor does the concept of
stable intersection as introduced in \cite{RGST} work for arbitrary
ambient spaces as can be seen in example \ref{example-rigidcurve}.
Therefore we introduce the notion of affine Cartier divisors on
tropical cycles as piecewise integer affine linear functions
modulo globally affine linear functions and define a bilinear
intersection product of Cartier divisors and cycles. We then prove
the commutativity of this product and a projection formula for
push-forwards of cycles and pull-backs of Cartier divisors. In the
second part (sections \ref{sec-cycles} -
\ref{sec-rationalequivalence}) we generalize the theory developed
in the first part to abstract cycles which are abstract polyhedral
complexes modulo refinements with affine cycles as local building
blocks. Again, abstract tropical varieties are just cycles with
non-negative weights. In both the affine and abstract case a remarkable
difference to the classical situation occurs: We can define the
mentioned intersection products on the level of cycles, i.e. we
can intersect Cartier divisors with cycles and obtain a
well-defined cycle --- not only a cycle class up to rational
equivalence as it is the case in classical algebraic geometry.
However, for simplifying the computations of concrete enumerative
numbers we introduce a notion of rational equivalence of cycles in
section \ref{sec-rationalequivalence}. In the third part (section
\ref{sec-stableintersection}) we finally use our theory to define
the intersection product of two cycles with ambient space $\RR^n$.
Here again it is remarkable that we can define these intersections
--- even for self-intersections --- on the level of cycles. We suppose
this intersection product to be identical with the \df{stable
intersection} discussed in \cite{M} and \cite{RGST} though we
could not prove it yet.

There are three more articles related to our work that we want to
mention: In \cite{K} the author studies the relations between the
intersection products of toric varieties and the tropical
intersection product on $\RR^n$ in the case of transversal
intersections. This article is closely related to \cite{FS}: In
this work the authors give a description of the Chow cohomology of
a complete toric variety  in terms of \df{Minkowski weights}.
These objects --- representing cocycles in the toric variety ---
are affine tropical cycles in $\RR^n$ according to our definition.
Moreover, there is an intersection product of these Minkowski
weights corresponding to the cup product of the associated
cocycles that can be calculated via a \df{fan displacement rule}.
This rule equals the stable intersection of tropical cycles in
$\RR^n$ mentioned above for the case of affine cycles. But there
are also discrepancies between these two interpretations of
Minkowski weights: Morphisms of toric varieties as well as
morphisms of affine tropical cycles are just given by integer
linear maps. However, the requirements for the fans are quite
different for both kinds of morphisms. Also the functorial
behavior is totally different for both interpretations: Regarded
as toric cocycles, Minkowski weights have pull-backs along
morphisms, whereas interpreted as affine tropical cycles they
admit push-forwards. In \cite{ST} the authors study homomorphisms
of tori and their induced morphisms of toric varieties and
tropical varieties, respectively. Generically finite morphisms in
this context are closely related to push-forwards of tropical
cycles as defined in construction \ref{constr-push}.

We would like to thank our advisor Andreas Gathmann for numerous
helpful discussions and his inspiring ideas that made this paper
possible.

  \section{Affine tropical cycles} \label{sec-affinecycles}

In this section we will briefly summarize the definitions and some
properties of our basic objects. We refer to \cite{GKM} for more
details (but note that we use a slightly more general definition of fan).\\
In the following sections $\Lambda$ will denote a finitely
generated free abelian group, i.e. a group isomorphic to $\ZZ^r$
for some $r \in \NN$, and $V := \Lambda \otimes_\ZZ \RR$ the
associated real vector space containing $\Lambda$ as a lattice. We
will denote the dual lattice in the dual vector space by
$\Lambda^{\vee} \subseteq V^{\vee}$.

\begin{definition}[Cones] \label{def-cones}
  A \df{cone} in $V$ is a subset $\sigma \subseteq V$ that can be
  described by finitely many linear integral equalities and
  inequalities, i.e. a set of the form
  $$\sigma = \left\{ x \in V | f_1(x)=0, \dots, f_r(x)=0, f_{r+1}(x) \geq 0, \dots, f_N(x) \geq 0
  \right\}$$
  for some linear forms $f_1, \dots, f_N \in \Lambda^{\vee}$. We
  denote by $V_\sigma$ the smallest linear subspace of $V$
  containing $\sigma$ and by $\Lambda_\sigma$ the lattice
  $V_\sigma \cap \Lambda$. We define the \df{dimension} of
  $\sigma$ to be the dimension of $V_\sigma$.
\end{definition}

\begin{definition}[Fans] \label{def-fans}
  A \df{fan} $X$ in $V$ is a finite set of
  cones in $V$ satisfying the following conditions:
  \begin{enumerate}
    \item \label{prop1} The intersection of any two cones in $X$ belongs to $X$ as well,
    \item \label{prop2} every cone $\sigma \in X$ is the disjoint union $\sigma =
    {\bigcup\limits^{\mbox{\Large .}}}_{\tau \in X: \tau \subseteq \sigma} \tau^{ri}$,
    where $\tau^{ri}$ denotes the relative interior of $\tau$, i.e. the interior of $\tau$ in $V_{\tau}$.
  \end{enumerate}
  We will denote the set of all $k$-dimensional cones of $X$ by
  $X^{(k)}$. The \df{dimension} of $X$ is defined to be the maximum of
  the dimensions of the cones in $X$. The fan $X$ is called
  \df{pure-dimensional} if each inclusion-maximal cone in $X$
  has this dimension. The union of all cones in $X$ will be denoted $|X|
  \subseteq V$. If $X$ is a fan of pure dimension $k$ then the
  cones $\sigma \in X^{(k)}$ are called \df{facets} of $X$. \medskip \\
  Let $X$ be a fan and $\sigma \in X$ a cone. A cone $\tau \in X$
  with $\tau \subseteq \sigma$ is called a \df{face} of $\sigma$. We
  write this as $\tau \leq \sigma$ (or $\tau < \sigma$ if in
  addition $\tau \subsetneq \sigma$ holds). Clearly we have
  $V_\tau \subseteq V_\sigma$ and $\Lambda_\tau \subseteq
  \Lambda_\sigma$ in this case. Note that $\tau < \sigma$ implies that
  $\tau$ is contained in a proper face (in the usual sense) of $\sigma$.
\end{definition}

\begin{construction}[Normal vectors] \label{constr-normalvectors}
  Let $\tau < \sigma$ be cones of some fan $X$ in $V$ with $\dim(\tau)=\dim(\sigma)-1$.
  This implies that there is a linear form $f \in \Lambda_{\sigma}^{\vee}$ that
  is zero on $\tau$, non-negative on $\sigma$ and not identically
  zero on $\sigma$. Let $u_{\sigma} \in \Lambda_{\sigma}$ be a vector
  generating $\Lambda_{\sigma} / \Lambda_{\tau} \cong \ZZ$ with
  $f(u_{\sigma}) > 0$.
  Note that its class $u_{\sigma / \tau} := \left[ u_{\sigma} \right] \in \Lambda_{\sigma} /
  \Lambda_{\tau}$ does not depend on the choice of
  $u_{\sigma}$. We call $u_{\sigma / \tau}$ the \df{(primitive)
  normal vector} of $\sigma$ relative to $\tau$.
\end{construction}

\begin{definition}[Subfans] \label{def-subfans}
  Let $X,Y$ be fans in $V$. $Y$ is called
  a \df{subfan} of $X$ if for every cone $\sigma \in Y$ there
  exists a cone $\sigma' \in X$ such that $\sigma \subseteq
  \sigma'$. In this case we write $Y \subfan X$ and
  define a map $C_{Y,X}:Y \rightarrow X$ that maps a cone
  $\sigma \in Y$ to the unique inclusion-minimal cone
  $\sigma' \in X$ with $\sigma \subseteq \sigma'$.
\end{definition}

\begin{definition}[Weighted fans] \label{def-wfans}
  A \df{weighted fan} $(X, \omega_X)$ of dimension $k$ in $V$ is a
  fan $X$ in $V$ of pure dimension $k$,
  together with a map $\omega_X : X^{(k)} \rightarrow \ZZ$.
  The number $\omega_X(\sigma)$ is called the \df{weight} of the facet $\sigma
  \in X^{(k)}$. For simplicity we usually write $\omega(\sigma)$
  instead of $\omega_X(\sigma)$. Moreover, we want to consider
  the \df{empty fan} $\emptyset$ to be a weighted fan of
  dimension $k$ for all $k$. Furthermore, by abuse of notation
  we simply write $X$ for the weighted fan $(X,\omega_X)$
  if the weight function $\omega_X$ is clear from the context.
\end{definition}

\begin{definition}[Tropical fans] \label{def-tropfans}
  A \df{tropical fan} of dimension $k$ in $V$ is a weighted fan
  $(X,\omega_X)$ of dimension $k$ satisfying the following \df{balancing condition}
  for every $\tau \in X^{(k-1)}$:
    $$\sum_{\sigma: \tau<\sigma} \omega_X(\sigma) \cdot u_{\sigma / \tau} =
    0 \in V / V_{\tau}.$$
  Let $(X,\omega_X)$ be a weighted fan of dimension $k$ in $V$ and $X^{*}$ the
  fan $$X^{*}:=\{ \tau \in X| \tau \leq \sigma \text { for some
  facet } \sigma \in X \text{ with } \omega_X(\sigma) \neq 0\}.$$
  $(X^{*},\omega_{X^{*}}):=(X^{*},\omega_X|_{(X^{*})^{(k)}})$ is called
  the \df{non-zero part} of $X$ and is again a weighted fan of dimension $k$ in $V$
  (note that $X^{*}=\emptyset$ is possible). Obviously $(X^{*},\omega_{X^{*}})$
  is a tropical fan if and only if $(X,\omega_X)$ is one. We call a weighted fan $(X,\omega_X)$
  \df{reduced} if all its facets have non-zero weight, i.e. if $(X,\omega_X)=(X^{*},\omega_{X^*})$ holds.
\end{definition}

\begin{remark} \label{rem-balancingcond}
  Let $(X,\omega_X)$ be a tropical fan of dimension $k$ and let
  $\tau \in X^{(k-1)}$. Let $\sigma_1, \dots, \sigma_N$ be all
  cones in $X$ with $\sigma_i > \tau$. For all $i$ let $v_{\sigma_i / \tau} \in \Lambda$
  be a representative of the primitive normal vector $u_{\sigma_i / \tau} {\in \Lambda / \Lambda_{\tau}}$.
  By the above balancing condition we have $\sum_{i=1}^{N}
  \omega_X(\sigma_i) \cdot v_{\sigma_i / \tau} = \lambda_{\tau}$ for some $\lambda_{\tau} \in \Lambda_{\tau}$.
  Obviously we have $\lambda_{\tau} = {\text{gcd}(\omega_X(\sigma_1), \dots,
  \omega_X(\sigma_N)) \cdot \widetilde{\lambda}_{\tau}}$ for some further $\widetilde{\lambda}_{\tau} \in \Lambda_{\tau}$.
  We can represent the greatest common divisor
  by a linear combination $\text{gcd}(\omega_X(\sigma_1), \dots,
  \omega_X(\sigma_N)) = \alpha_1 \omega_X(\sigma_1) + \dots + \alpha_N
  \omega_X(\sigma_N)$ with $\alpha_1, \dots, \alpha_N \in \ZZ$ and define $${\widetilde{v}_{\sigma_i / \tau}:=v_{\sigma_i /
  \tau}- \alpha_i \cdot \widetilde{\lambda}_{\tau}}$$ for all $i$.
  Note that $\widetilde{v}_{\sigma_i / \tau}$ is a representative of $u_{\sigma_i /
  \tau}$, too. Replacing all $v_{\sigma_i / \tau}$ by $\widetilde{v}_{\sigma_i /
  \tau}$ we can always assume that $\sum_{i=1}^{N} \omega_X(\sigma) \cdot v_{\sigma /
  \tau}=0 \in \Lambda$.
\end{remark}

\begin{definition}[Refinements] \label{def-refinements}
  Let $(X,\omega_X)$ and $(Y,\omega_Y)$ be weighted fans in $V$. We
  call $(Y,\omega_Y)$ a \df{refinement} of $(X,\omega_X)$ if the
  following holds:
  \begin{enumerate}
    \item $Y^{*} \subfan X^{*}$,
    \item \label{prop2b} $|Y^{*}|=|X^{*}|$ and
    \item $\omega_Y(\sigma)=\omega_X(C_{Y^{*},X^{*}}(\sigma))$ for
    every $\sigma \in (Y^{*})^{(\dim(Y))}$.
  \end{enumerate}
  Note that property \ref{prop2b} implies that either $X^{*}=Y^{*}=\emptyset$
  or $\dim(X)=\dim(Y)$. We call two weighted fans $(X,\omega_X)$ and
  $(Y,\omega_Y)$ in $V$ \df{equivalent} (write
  $(X,\omega_X) \sim (Y,\omega_Y)$) if they have a common
  refinement. Note that $(X,\omega_X)$and $(X^{*},\omega_X|_{(X^{*})^{(\dim(X))}})$ are
  always equivalent.
\end{definition}

\begin{remark} \label{rem-refinement}
  Note that for a weighted fan $(X,\omega_X)$ of dimension $k$ and a
  refinement $(Y,\omega_Y)$ we have the following
  two properties:
  \begin{enumerate}
    \item $|X^{*}|=|Y^{*}|$, i.e. the support $|X^{*}|$ is well-defined
    on the equivalence class of $X$,
    \item for every cone $\tau \in Y^{(k-1)}$ there are exactly two cases
    that can occur: Either $\dim C_{Y,X}(\tau)=k$ or $\dim C_{Y,X}(\tau)=k-1$.
    In the first case all cones $\sigma \in Y^{(k)}$ with
    $\sigma>\tau$
    must be contained in $C_{Y,X}(\tau)$. Thus there are precisely two such cones
    $\sigma_1$ and $\sigma_2$ with $\omega_Y(\sigma_1)=\omega_Y(\sigma_2)$ and
    $u_{\sigma_1 / \tau} = - u_{\sigma_2 / \tau}$. In the second
    case we have a 1:1 correspondence between cones $\sigma \in
    Y^{(k)}$ with $\tau<\sigma$ and cones $\sigma' \in X^{(k)}$
    with $C_{Y,X}(\tau) < \sigma'$ preserving weights and normal vectors.
  \end{enumerate}
\end{remark}

\begin{construction}[Refinements] \label{constr-refinement}
  Let $(X,\omega_X)$ be a weighted fan and $Y$ be any fan in $V$ with
  $|X| \subseteq |Y|$. Let $P:= \{ \sigma \cap \sigma'
  | \sigma \in X, \sigma' \in Y \}$. In general $P$ is not a fan in $V$
  as can be seen in the following example:\\
  \medskip\\
  \begin{minipage}{\linewidth}
  \begin{center}
  \input{pic/Ex_XcapY.pstex_t}\\
  Fans $X$ and $Y$ such that $\{ \sigma \cap \sigma'| \sigma \in X, \sigma' \in Y \}$ is not a fan.
  \end{center}
  \end{minipage}
  \medskip\\
  Here $P$ contains $\tau_1'=\sigma_2 \cap \sigma_1'$, but also
  $\tau_2=\sigma_1 \cap \sigma_2'$ and $\tau_3=\sigma_3 \cap
  \sigma_2'$. Hence property (b) of definition \ref{def-fans}
  is not fulfilled. To resolve this, we define $$X \cap Y :=\{ \sigma \in P| \nexists \enspace \tau \in
  P^{(\dim(\sigma))} \text{ with } \tau \subsetneq \sigma \}.$$
  Note that $X \cap Y$ is now a fan in $V$.
  We can make it into a weighted fan by setting $\omega_{X \cap Y}(\sigma):=
  \omega_X(C_{X \cap Y,X}(\sigma))$ for all $\sigma \in (X \cap Y)^{(\dim(X))}$.
  Then $(X \cap Y,\omega_{X \cap Y})$ is a refinement of $(X,\omega_X)$.
  Note that if $(X,\omega_X)$ and $(Y,\omega_Y)$ are both weighted fans
  and $|X|=|Y|$ we can form both intersections $X \cap Y$ and $Y \cap X$.
  Of course, the underlying fans are the same in both cases, but the
  weights may differ since they are always induced by the first complex.
  \medskip \\
  The following setting is a special case of this construction:
  Let $(X,\omega_X)$ be a weighted fan of dimension $k$ in $V$
  and let $f \in \Lambda^{\vee}$ be a non-zero linear form.
  Then we can construct a refinement of $(X,\omega_X)$ as follows:
  $$H_{f} := \left\{ \{ x \in V | f(x) \leq 0 \}, \{ x \in V | f(x) = 0 \},
  \{ x \in V | f(x) \geq 0 \} \right\}$$ is a fan in $V$
  with $|H_{f}| = V$. Thus we have $|X| \subseteq |H_{f}|$
  and by our above construction we get a refinement
  $(X_{f},\omega_{X_{f}}) := (X \cap H_{f}, \omega_{X \cap H_{f}})$ of $X$.
\end{construction}

Obviously we still have to answer the question if the equivalence
of weighted fans is indeed an equivalence relation and if this
notion of equivalence is well-defined on tropical fans. We will do
this in the following lemma:

\begin{lemma} \label{lemma-equivrel}
  \begin{enumerate}
    \item The relation ``$\sim$'' is an equivalence relation
    on the set of $k$-dimensional weighted fans in $V$.
    \item If $(X,\omega_X)$ is a weighted fan of dimension $k$ and
    $(Y,\omega_Y)$ is a refinement then $(X,\omega_X)$ is a
    tropical fan if and only if $(Y,\omega_Y)$ is one.
  \end{enumerate}
\end{lemma}
\begin{proof}
  Recall that a fan and its non-zero part are always equivalent
  and that a weighted fan $X$ is tropical if and only if its non-zero part
  $X^{*}$ is. Thus we may assume that all our fans are
  reduced and the proof is the same as in \cite[section 2]{GKM}.
\end{proof}

Having done all these preparations we are now able to introduce
the most important objects for the succeeding sections:

\begin{definition}[Affine cycles and affine tropical varieties] \label{def-affcyclesvarieties}
  Let $(X,\omega_X)$ be a tropical fan of dimension $k$ in $V$. We denote by $[(X,\omega_X)]$ its equivalence class
  under the equivalence relation ``$\sim$'' and by $Z_k^{\aff}(V)$ the set of equivalence classes
    $$Z_k^{\aff}(V):=\{ [(X,\omega_X)] | (X,\omega_X) \text{ tropical fan of dimension } k \text{ in } V\}.$$
  The elements of $Z_k^{\aff}(V)$ are called \df{affine (tropical) $k$-cycles} in $V$. A $k$-cycle
  $[(X,\omega_X)]$ is called an \df{affine tropical variety} if
  $\omega_X(\sigma) \geq 0$ for every $\sigma \in X^{(k)}$.
  Note that the last property is independent of the choice of the
  representative of $[(X,\omega_X)]$. Moreover, note that
  $0:=[\emptyset] \in Z_k^{\aff}(V)$ for every $k$. We define
  $|[(X,\omega_X)]|:=|X^{*}|$. This definition is well-defined by
  remark \ref{rem-refinement}.
\end{definition}

\begin{construction}[Sums of affine cycles] \label{constr-sum}
  Let $[(X,\omega_X)]$ and $[(Y,\omega_Y)]$ be $k$-cycles in $V$.
  We would like to form a fan $X+Y$ by taking the union $X \cup
  Y$, but obviously this collection of cones is in general not a fan.
  Using appropriate refinements we can resolve this problem:
  Let ${f_1(x) \geq 0, \dots, f_{N_1}(x) \geq
  0}$, ${f_{N_1+1}(x) = 0, \dots, f_N(x)=0}$ and ${g_1(x) \geq 0, \dots, g_{M_1}(x) \geq
  0}$, ${g_{M_1+1}(x) = 0, \dots, g_M(x)=0}$ be all different equalities
  and inequalities occurring in the descriptions of all the cones belonging to $X$
  and $Y$ respectively. Using construction \ref{constr-refinement}
  we get refinements $$\widetilde{X}:=X \cap H_{f_1} \cap \dots \cap H_{f_N} \cap H_{g_1} \cap \dots \cap H_{g_M}$$
  of $X$ and $$\widetilde{Y}:=Y \cap H_{f_1} \cap \dots \cap H_{f_N} \cap H_{g_1} \cap \dots \cap H_{g_M}$$ of
  $Y$ (note that the final refinements do not depend on the order of the single refinements).
  A cone occurring in $\widetilde{X}$ or $\widetilde{Y}$ is then of the form
  $$ \sigma = \left\{ \left. \begin{array}{ccc}
      {f_i(x) \leq 0}, & {f_{j}(x) = 0}, & {f_k(x) \geq 0}, \\
      {g_{i'}(x) \leq 0}, & {g_{j'}(x) = 0}, & {g_{k'}(x) \geq 0} \end{array} \right|
      \begin{array}{ccc}
      i \in I, &  j \in J, & k \in K, \\i' \in I', & j' \in J', & k' \in K'
      \end{array}
    \right\}$$
  for some partitions $I \cup J \cup K = \{1,\dots,N\}$ and $I' \cup J' \cup K' =
  \{1,\dots,M\}$. Now, all these cones $\sigma$ belong
  to the fan $H_{f_1} \cap \dots \cap H_{f_N} \cap H_{g_1} \cap \dots \cap
  H_{g_M}$ as well and hence $\widetilde{X} \cup \widetilde{Y}$
  fulfills definition \ref{def-fans}.
  Thus, now we can define the \df{sum of $X$ and $Y$} to be $X+Y := \widetilde{X} \cup \widetilde{Y}$
  together with weights $\omega_{X+Y} (\sigma) := \omega_{\widetilde{X}} (\sigma)+\omega_{\widetilde{Y}}
  (\sigma)$ for every facet of $X+Y$ (we set $\omega_{\square} (\sigma):=0$ if $\sigma$ does not occur
  in $\square \in \{ \widetilde{X}, \widetilde{Y} \} )$. By
  construction, $(X+Y,\omega_{X+Y})$ is again a tropical fan of
  dimension $k$. Moreover, enlarging the sets $\{f_i\}, \{g_j\}$ by more
  (in)equalities just corresponds to refinements of $X$ and $Y$
  and only leads to a refinement of $X+Y$. Thus, replacing the set
  of relations by another one that also describes the cones in
  $X$ and $Y$, or replacing $X$ or $Y$ by refinements keeps the
  equivalence class $[(X+Y,\omega_{X+Y})]$ unchanged, i.e. taking
  sums is a well-defined operation on cycles.
\end{construction}

This construction immediately leads to the following lemma:

\begin{lemma} \label{lemma-group}
  $Z_k^{\aff}(V)$ together with the operation ``+'' from construction \ref{constr-sum} forms an abelian group.
\end{lemma}
\begin{proof}
  The class of the empty fan $0=[\emptyset]$ is the neutral
  element of this operation and $[(X,-\omega_X)]$ is the inverse element of
  $[(X,\omega_X)] \in Z_k^{\aff}(V)$.
\end{proof}

Of course we do not want to restrict ourselves to cycles situated
in some $\RR^n$. Therefore we give the following generalization of
definition \ref{def-affcyclesvarieties}:

\begin{definition} \label{def-gencyclegroup}
  Let $X$ be a fan in $V$. An \df{affine $k$-cycle in
  $X$} is an element $[(Y,\omega_Y)]$ of $Z_k^{\aff}(V)$
  such that $|Y^{*}| \subseteq |X|$. We denote by
  $Z_k^{\aff}(X)$ the set of $k$-cycles in $X$.
  Note that $\left( Z_k^{\aff}(X),+ \right)$ is a subgroup of
  $(Z_k^{\aff}(V),+)$. The elements of the group $Z_{\dim X -1}^{\aff}(X)$ are called
  \df{Weil divisors} on $X$.
  If $[(X,\omega_X)]$ is a cycle in $V$ then
  $Z_k^{\aff}([(X,\omega_X)]) := Z_k^{\aff}(X^{*})$.
\end{definition}

  \section{Affine Cartier divisors and their associated Weil divisors} \label{sec-affinecartierdivisors}

\begin{definition}[Rational functions] \label{def-rationalfunctions}
  Let $C$ be an affine $k$-cycle. A \df{(non-zero) rational function on $C$} is a continuous piecewise linear function $\varphi : |C| \rightarrow \RR$, i.e. there exists a representative $(X,\omega_X)$ of $C$ such that on each cone $\sigma \in X$, $\varphi$ is the restriction of an integer affine linear function $\varphi|_\sigma = \lambda + c$, $\lambda \in \Lambda^\vee_\sigma$, $c \in \RR$. Obviously, $c$ is the same on all faces by $c = \varphi(0)$ and $\lambda$ is uniquely determined by $\varphi$ and therefore denoted by $\varphi_\sigma := \lambda$. \\
  The \df{set of (non-zero) rational functions of $C$} is denoted by $\calK^*(C)$.
\end{definition}

\begin{remark}[The zero function and restrictions to subcycles] \label{remark-zerofunction}
  The ``zero'' function can be thought of being the constant function $-\infty$, therefore $\calK(C) := \calK^*(C) \cup \{-\infty\}$. With respect to the operations $\max$ and $+$, $\calK(C)$ is a semifield. \\
  Let us note an important difference to the classical case: Let $D$ be an arbitrary subcycle of $C$ and $\varphi \in \calK^*(C)$. Then $\varphi |_{|D|} \in \calK^*(D)$, whereas in the classical case it might become zero. This will be crucial for defining intersection products not only modulo rational equivalence. On the other hand, the definition of rational functions given above, requiring the function to be defined everywhere, seems to be restrictive when compared to the classical case, even so ``being defined'' does not imply ``being regular'' tropically. In some cases (see remark \ref{remark-pushforwardratequiv}) it would be desirable to generalize our definition while preserving the above restriction property.
\end{remark}

As in the classical case, each non-zero rational function
$\varphi$ on $C$ defines a Weil divisor, i.e. a cycle in $Z_{\dim
C - 1}^{\aff}(C)$. The idea of course should be to describe the
``zeros'' and ``poles'' of $\varphi$. A naive approach could be to
consider the graph of $\varphi$ in $V \times \RR$ and ``intersect
it with $V \times \{-\infty\}$ and $V \times \{+\infty\}$''.
However, our function $\varphi$ takes values only in $\RR$, in
fact. On the other hand, the graph of $\varphi$ is not a tropical
object as it is not balanced: Where $\varphi$ is not linear, our
graph gets edges that might violate the balancing condition. So,
we first make the graph balanced by adding new faces in the
additional direction $(0, -1) \in V \times \RR$ and then apply our
naive approach. Let us make this precise.

\begin{construction}[The associated Weil divisor] \label{constr-associatedweildivisor}
  Let $C$ be an affine $k$-cycle in $V = \Lambda \otimes \RR$ and $\varphi \in \calK^*(C)$ a rational function on $C$. Let furthermore $(X, \omega)$ be a representative of $C$ on whose faces $\varphi$ is affine linear. Therefore, for each cone $\sigma \in X$, we get a cone $\tilde{\sigma} := (\id \times \varphi_\sigma)(\sigma)$ in $V \times \RR$ of the same dimension. Obviously, $\Gamma_\varphi := \{\tilde{\sigma} | \sigma \in X\}$ forms a fan which we can make into a weighted fan $(\Gamma_\varphi, \tilde{\omega})$ by $\tilde{\omega}(\tilde{\sigma}) := \omega(\sigma)$. Its support is just the set-theoretic graph of $\varphi - \varphi(0)$ in $|X| \times \RR$.\\
  For $\tau < \sigma$ with $\dim(\tau)=\dim(\sigma)-1$ let $v_{\sigma / \tau} \in \Lambda$ be a representative of the normal vector $u_{\sigma / \tau}$. Then, $\left(v_{\sigma / \tau}, \varphi_\sigma(v_{\sigma / \tau})\right) \in \Lambda \times \ZZ$ is a representative of the normal vector $u_{\tilde{\sigma} / \tilde{\tau}}$. Therefore, summing around a cone $\tilde{\tau}$ with $\dim \tilde{\tau} = \dim \tau = k - 1$, we get
  $$
    \sum_{\twolines{\tilde{\sigma} \in \Gamma_\varphi^{(k)}}{\tilde{\tau} < \tilde{\sigma}}}
        \tilde{\omega}(\tilde{\sigma}) \left(v_{\sigma / \tau}, \varphi_\sigma(v_{\sigma / \tau}) \right)
        =
        \left(\sum_{\twolines{\sigma \in X^{(k)}}{\tau < \sigma}} \omega(\sigma) v_{\sigma / \tau},
        \sum_{\twolines{\sigma \in X^{(k)}}{\tau < \sigma}} \varphi_\sigma(\omega(\sigma) v_{\sigma / \tau}) \right).
  $$
  From the balancing condition for $(X, \omega)$ it follows that $\sum_{\sigma \in X^{(k)} : \tau < \sigma} \omega(\sigma) v_{\sigma / \tau} \in V_\tau$, which also means $\left( \sum_{\sigma \in X^{(k)} : \tau < \sigma} \omega(\sigma) v_{\sigma / \tau}, \varphi_\tau \left( \sum_{\sigma \in X^{(k)} : \tau < \sigma} \omega(\sigma) v_{\sigma / \tau} \right) \right) \in V_{\tilde{\tau}}$. Therefore, modulo $V_{\tilde{\tau}}$, our first sum equals
  $$
    \left( 0,
    \sum_{\twolines{\sigma \in X^{(k)}}{\tau < \sigma}} \varphi_\sigma(\omega(\sigma) v_{\sigma / \tau}) -
    \varphi_\tau\Big( \sum_{\twolines{\sigma \in X^{(k)}}{\tau < \sigma}} \omega(\sigma) v_{\sigma / \tau}\Big) \right) \in V \times \RR.
  $$
  So, in order to ``make $(\Gamma_\varphi, \tilde{\omega})$ balanced at $\tilde{\tau}$'', we add the cone $\vartheta := \tilde{\tau} + (\{0\} \times \RR_{\leq 0})$ with weight $\tilde{\omega}(\vartheta) = \sum_{\sigma \in X^{(k)} : \tau < \sigma} \varphi_\sigma(\omega(\sigma) v_{\sigma / \tau}) -
    \varphi_\tau\Big( \sum_{\sigma \in X^{(k)} : \tau < \sigma} \omega(\sigma) v_{\sigma / \tau}\Big) $. As obviously $[(0, -1)] = u_{\vartheta / \tilde{\tau}} \in (V \times \RR)/V_{\tilde{\tau}}$, the above calculation shows that then the balancing condition around $\tilde{\tau}$ holds. In other words, we build the new fan $(\Gamma'_\varphi, \tilde{\omega}')$, where
  \begin{eqnarray*}
    \Gamma'_\varphi & := & \Gamma_\varphi \cup \left\{ \tilde{\tau} + (\{0\} \times \RR_{\leq 0}) | \tilde{\tau} \in \Gamma_\varphi \setminus \Gamma_{\varphi}^{(k)} \right\}, \\
    \tilde{\omega}'|_{\Gamma_{\varphi}^{(k)}} & := & \tilde{\omega}, \\
    \tilde{\omega}'(\tilde{\tau} + (\{0\} \times \RR_{\leq 0})) & := & \sum_{\twolines{\sigma \in X^{(k)}}{\tau < \sigma}} \varphi_\sigma(\omega(\sigma) v_{\sigma / \tau}) -
    \varphi_\tau\Big( \sum_{\twolines{\sigma \in X^{(k)}}{\tau < \sigma}} \omega(\sigma) v_{\sigma / \tau}\Big) \\ & & \text{ if } \dim \tilde{\tau} = k - 1.
  \end{eqnarray*}
  This fan is balanced around all $\tilde{\tau} \in \Gamma_{\varphi}^{(k-1)}$. We will show that it is also balanced at all ``new'' cones of dimension $k-1$ in proposition \ref{prop-balancingconditionandcommutativity}. \\
  We now think of intersecting this new fan with $V \times \{-\infty\}$ to get our desired Weil divisor (As our weights are allowed to be negative, we can forget about intersecting also with $V \times \{+\infty\}$). This construction leads to the following definition.
  \begin{figure}
    \input{pic/Ex_Graph.pstex_t}\\
    Construction of a Weil divisor.
  \end{figure}
\end{construction}

\begin{definition}[Associated Weil divisors] \label{def-associatedweildivisor}
    Let $C$ be an affine $k$-cycle in $V = \Lambda \otimes \RR$ and $\varphi \in \calK^*(C)$ a rational function on $C$. Let furthermore $(X, \omega)$ be a representative of $C$ on whose cones $\varphi$ is affine linear. We define $\divisor(\varphi) := \varphi \cdot C := [(\bigcup_{i=0}^{k-1} X^{(i)}, \omega_\varphi)] \in Z_{k-1}^{\aff}(C)$, where
    \begin{eqnarray*}
        \omega_\varphi : X^{(k-1)} & \rightarrow & \ZZ, \\
        \tau                                             & \mapsto     & \sum_{\twolines{\sigma \in X^{(k)}}{\tau < \sigma}} \varphi_\sigma(\omega(\sigma) v_{\sigma / \tau}) -
    \varphi_\tau\Big( \sum_{\twolines{\sigma \in X^{(k)}}{\tau < \sigma}} \omega(\sigma) v_{\sigma / \tau}\Big)
    \end{eqnarray*}
    and the $v_{\sigma / \tau}$ are arbitrary representatives of the normal vectors $u_{\sigma / \tau}$. \\
    Let $D$ be an arbitrary subcycle of $C$. By remark \ref{remark-zerofunction}, we can define $\varphi \cdot D:= \varphi|_{|D|} \cdot D$.
\end{definition}

\begin{remark}  \label{rem-welldefined}
    Obviously, $\omega_\varphi(\tau)$ is independent of the choice of the $v_{\sigma / \tau}$, as a different choice only differs by elements in $V_\tau$. \\
    Our definition does also not depend on the choice of a representative $(X, \omega)$: Let $(Y, \upsilon)$ be a refinement of $(X, \omega)$. For $\tau \in Y^{(k-1)}$, two cases can occur (see also remark \ref{rem-refinement}): Let $\tau' := C_{Y,X}(\tau)$. If $\dim \tau' = k$, there are precisely two cones at $\tau < \sigma_1, \sigma_2 \in Y^{(k)}$, which then fulfill $C_{Y,X}(\sigma_1) = C_{Y,X}(\sigma_2)$ and therefore $u_{\sigma_1 / \tau} = - u_{\sigma_2 / \tau}$, $\upsilon(\sigma_1) =\upsilon(\sigma_2)$ and $\varphi_{\sigma_1} = \varphi_{\sigma_2}$. It follows that $\upsilon_\varphi(\tau) = 0$. If $\dim \tau' = k-1$, $C_{Y,X}$ gives a one-to-one correspondence between $\{\sigma \in Y^{(k)} | \tau < \sigma\}$ and $\{\sigma' \in X^{(k)} | \tau' < \sigma'\}$ respecting weights and normal vectors, and we have $\varphi_{\sigma} = \varphi_{C_{Y,X}(\sigma)}$. It follows that $\upsilon_\varphi(\tau) = \omega_\varphi(\tau')$. So the two weighted fans we obtain are equivalent.
\end{remark}

\begin{remark}[Affine linear functions and sums] \label{rem-affinelinearfunctionsandsums}
  Let $\varphi \in \calK^*(C)$ be globally affine linear, i.e. $\varphi = \lambda|_{|C|} + c$ for some $\lambda \in \Lambda^\vee$, $c \in \RR$. Then obviously $\varphi \cdot C = 0$. \\
  Let furthermore $\psi \in \calK^*(C)$ be another rational function on $C$. From $\varphi_\sigma + \psi_\sigma = (\varphi + \psi)_\sigma$ it follows that $(\varphi + \psi) \cdot C = \varphi \cdot C + \psi \cdot C$.
\end{remark}

\pagebreak
\begin{proposition}[Balancing Condition and Commutativity] \label{prop-balancingconditionandcommutativity}
    \begin{enumerate}
    \item
    Let $C$ be an affine $k$-cycle in $V = \Lambda \otimes \RR$ and $\varphi \in \calK^*(C)$ a rational function on $C$. Then $\divisor(\varphi) = \varphi \cdot C$ is an equivalence class of tropical fans, i.e. its representatives are balanced.
    \item
    Let $\psi \in \calK^*(C)$ be another rational function on $C$. Then it holds $\psi \cdot (\varphi \cdot C) = \varphi \cdot (\psi \cdot C)$.
  \end{enumerate}
\end{proposition}

\begin{proof}
  (a): Let $(X, \omega)$ be a representative of $C$ on whose cones
  $\varphi$ is affine linear. Pick a $\theta \in X^{(k-2)}$. We choose
  an element $\lambda \in \Lambda^\vee$ with $\lambda|_{V_\theta} =
  \varphi_\theta$. By remark \ref{rem-affinelinearfunctionsandsums},
  we can go on with $\varphi - \lambda - \varphi(0) \in \calK^*(C)$
  instead of $\varphi$. By dividing out $V_\theta$, we can restrict
  ourselves to the situation $\dim X = 2$, $\theta = \{0\}$. \\
  By a further refinement (i.e. by cutting an possibly occuring
  halfspace into two pieces along an additional ray), we can assume
  that all cones $\sigma \in X$ are simplicial. Therefore each
  two-dimensional cone $\sigma \in X^{(2)}$ is generated by two
  unique rays $\tau, \tau' \in X^{(1)}$, i.e. $\sigma = \tau + \tau'$.
  We denote
  $$
    \chi(\sigma) := [\Lambda_\sigma : \Lambda_\tau + \Lambda_{\tau'}]
        = [\Lambda_\sigma : \ZZ u_{\tau / \{0\}} + \ZZ u_{\tau' / \{0\}}],
  $$
  where $u_{\tau / \{0\}}$ and $u_{\tau' / \{0\}}$ denote the
  primitive normal vectors introduced in construction
  \ref{constr-normalvectors}. Then we get
  $$
    [u_{\tau' / \{0\}}] = \chi(\sigma) u_{\sigma / \tau} \mod V_\tau.
  $$
  This equation can be shown for example as follows: The linear
  extension of the following function
  \begin{eqnarray*}
    \text{index} : \Lambda_\sigma \setminus \Lambda_\tau & \rightarrow & \ZZ, \\
                   v                                     & \mapsto     &
                     [\Lambda_\sigma : \ZZ u_{\tau / \{0\}} + \ZZ v]
  \end{eqnarray*}
  to $\Lambda_\sigma$ is in fact trivial on $\Lambda_\tau$.
  Therefore it can also be considered as a function on
  $\Lambda_\sigma / \Lambda_\tau$. But by definitions
  we know $\text{index}(u_{\sigma / \tau}) = 1$ (as $u_{\tau / \{0\}}$
  and any representative of $u_{\sigma / \tau}$ form a lattice basis
  of $\Lambda_\sigma$) and $\text{index}(u_{\tau' / \{0\}}) =
  \chi(\sigma)$, which proves the claim. \\
  This means that we can rewrite the balancing condition of $X$
  around $\tau \in X^{(1)}$ only using the vectors generating
  the rays, namely
  \begin{eqnarray*}
    \sum_{\twolines{\tau' \in X^{(1)}}{\tau + \tau' \in X^{(2)}}}
      \frac{\omega(\sigma)}{\chi(\sigma)} u_{\tau' / \{0\}} & \in & V_\tau \\
    & = & \lambda_\tau u_{\tau / \{0\}},
  \end{eqnarray*}
  where $\lambda_\tau$ is a coefficient in $\RR$ and $\sigma$
  denotes $\tau + \tau'$ in such sums. Of course, we can also
  compute the weight $\omega_\varphi(\tau)$ of $\tau$ in
  $\divisor(\varphi)$:
  $$
    \omega_\varphi(\tau) =
      \left(\sum_{\twolines{\tau' \in X^{(1)}}{\tau + \tau' \in X^{(2)}}}
      \frac{\omega(\sigma)}{\chi(\sigma)} \varphi(u_{\tau' / \{0\}})\right)
      - \lambda_\tau \varphi(u_{\tau / \{0\}})
  $$
  Let us now check the balancing condition of $\varphi \cdot C$
  around $\{0\}$ by plugging in these equations. We get
  \begin{eqnarray*}
    \sum_{\tau \in X^{(1)}} \omega_\varphi(\tau) u_{\tau / \{0\}}
      & = &
      \sum_{\twolines{\tau, \tau' \in X^{(1)}}{\tau + \tau' \in X^{(2)}}}
      \frac{\omega(\sigma)}{\chi(\sigma)} \varphi(u_{\tau' / \{0\}}) u_{\tau / \{0\}} \\
    & & -
      \sum_{\tau \in X^{(1)}} \lambda_\tau \varphi(u_{\tau / \{0\}}) u_{\tau / \{0\}}.
  \end{eqnarray*}
  By Commuting $\tau$ and $\tau'$ in the first summand we get
  \begin{eqnarray*}
    \sum_{\tau \in X^{(1)}} \omega_\varphi(\tau) u_{\tau / \{0\}}
      & = &
      \sum_{\twolines{\tau, \tau' \in X^{(1)}}{\tau + \tau' \in X^{(2)}}}
      \frac{\omega(\sigma)}{\chi(\sigma)} \varphi(u_{\tau / \{0\}}) u_{\tau' / \{0\}} \\
    & & -
      \sum_{\tau \in X^{(1)}}
      \lambda_{\tau} \varphi(u_{\tau / \{0\}}) u_{\tau / \{0\}} \\
    & = &
      \sum_{\tau \in X^{(1)}} \varphi(u_{\tau / \{0\}})
      \underbrace{\left(\left(
      \sum_{\twolines{\tau' \in X^{(1)}}{\tau + \tau' \in X^{(2)}}}
      \frac{\omega(\sigma)}{\chi(\sigma)} u_{\tau' / \{0\}}\right)
      - \lambda_{\tau} u_{\tau / \{0\}}\right)
      }_{= 0 \text{ (balancing condition around $\tau$)}} \\
    & = & 0.
  \end{eqnarray*}
  This finishes the proof of (a). \\
  (b): Let $(X, \omega)$ be a representative of $C$ on whose cones
  $\varphi$ and $\psi$ are affine linear. Pick a $\theta \in
  X^{(k-2)}$. By the same reduction steps as in case (a), we can
  again restrict ourselves to $\dim X = 2$, $\theta = \{0\}$. With
  the notations and trick as in (a) we get
  $$
    \omega_{\varphi, \psi} (\{0\})
      = \sum_{\twolines{\tau, \tau' \in X^{(1)}}{\tau + \tau' \in X^{(2)}}}
      \frac{\omega(\sigma)}{\chi(\sigma)}
      \varphi(u_{\tau' / \{0\}}) \psi(u_{\tau / \{0\}})
      = \omega_{\psi, \varphi} (\{0\}),
  $$
  which finishes part (b).
\end{proof}

\begin{definition}[Affine Cartier divisors] \label{def-affinecartierdivisors}
    Let $C$ be an affine $k$-cycle. The subgroup of globally affine linear functions in $\calK^*(C)$ with respect to $+$ is denoted by $\calO^*(C)$. We define the \df{group of affine Cartier divisors of $C$} to be the quotient group $\Div(C) := \calK^*(C) / \calO^*(C)$. \\
    Let $[\varphi] \in \Div(C)$ be a Cartier divisor. By remark \ref{rem-affinelinearfunctionsandsums}, the associated Weil divisor \linebreak $\divisor([\varphi]) := \divisor(\varphi)$ is well-defined. We therefore get a bilinear mapping
    \begin{eqnarray*}
        \cdot \; : \Div(C) \times Z_k^{\aff}(C) & \rightarrow & Z_{k-1}^{\aff}(C), \\
        ([\varphi] , D)      & \mapsto & [\varphi] \cdot D = \varphi \cdot D,
    \end{eqnarray*}
    called \df{affine intersection product}.
\end{definition}

\begin{example}[Self-intersection of hyperplanes] \label{example-selfintersectionofhyperplanes}
    Let $\Lambda = \ZZ^n$ (and thus $V = \RR^n$), let $e_1, \ldots, e_n$ be the standard basis vectors in $\ZZ^n$ and $e_0 := - e_1 - \cdots - e_n$. By abuse of notation our ambient cycle is $\RR^n := [(\{\RR^n\}, \omega(\RR^n)= 1)]$. Let us consider the ``linear tropical polynomial'' $h = x_1 \oplus \cdots \oplus x_n \oplus 0 = \max\{x_1, \ldots, x_n, 0\} : \RR^n \rightarrow \RR$. Obviously, $h$ is a rational function in the sense of definition \ref{def-rationalfunctions}: For each subset $I \subsetneq \{0, 1, \ldots, n\}$ we denote by $\sigma_I$ the simplicial cone of dimension $|I|$ generated by the vectors $-e_i$ for $i \in I$. Then $h$ is integer linear on all $\sigma_I$, namely
    $$
        h|_{\sigma_I}(x_1, \ldots, x_n) =
        \left\{ \begin{array}{ll}
            0   &   \text{if } 0 \notin I, \\
            x_i &  \text{if there exists an } i \in \{1, \ldots, n\} \setminus I.
        \end{array} \right.
    $$
    Let $L^n_k$ be the $k$-dimensional fan consisting of all cones $\sigma_I$ with $|I| \leq k$ and weighted with the trivial weight function $\omega_{L^n_k}$. Then $L^n_n$ is a representative of $\RR^n$ fulfilling the conditions of definition \ref{def-rationalfunctions}. We want to show
    \[
        \underbrace{h \cdot \cdots \cdot h}_{k \text{ times}} \cdot \RR^n = [L^n_{n-k}]. \tag{$*$}
    \]
    This follows inductively from $h \cdot [L^n_{k+1}] = [L^n_k]$, so it remains to compute $\omega_{L^n_{k+1}, h}(\sigma_I)$ for all $I$ with $|I| = k < n$. Let $J := \{0, 1, \ldots, n\} \setminus I$. Obviously, the $(k+1)$-dimensional cones of $L^n_{k+1}$ containing $\sigma_I$ are precisely the cones $\sigma_{I \cup \{j\}}, j \in J$. Moreover, $-e_j$ is a representative of the normal vector $u_{\sigma_{I \cup \{j\}} / \sigma_I}$. Note also that for all $i \in I', I' \subsetneq \{0, 1, \ldots, n\}$ we have $h_{\sigma_{I'}}(-e_i) = h|_{\sigma_{I'}}(-e_i) = h(-e_i)$. Using this we compute
    \begin{eqnarray*}
        \omega_{L^n_{k+1}, h}(\sigma_I)
            & = &
            \sum_{j \in J} \underbrace{\omega_{L^n_{k+1}}(\sigma_{I \cup \{j\}})}_{=1} h_{\sigma_{I \cup \{j\}}}(-e_j) \\
        & &
            \text{} + h_{\sigma_I} \bigg(\underbrace{\sum_{j \in J} \underbrace{\omega_{L^n_{k+1}}(\sigma_{I \cup \{j\}})}_{=1} e_j}_{= -\sum_{i \in I} e_i}\bigg) \\
        & = &
            \sum_{j \in J} h(-e_j) + \sum_{i \in I} h(-e_i) \\
        & = &
            h(-e_0) + h(-e_1) + \cdots + h(-e_n) \\
        & = &
            1 + 0 + \cdots + 0 = 1 = \omega_{L^n_k} (\sigma_I),
    \end{eqnarray*}
    which implies $h \cdot [L^n_{k+1}] = [L^n_k]$ and also equation $(*)$. \\
    We can summarize this example as follows: Firstly, for a tropical polynomial $f$, the associated Weil divisor $f \cdot \RR^n$ coincides with the locus of non-differentiability $\calT(f)$ of $f$ (see \cite[section 3]{RGST}), and secondly, ``the $k$-fold self-intersection of a tropical hyperplane in $\RR^n$'' is given by its $(n-k)$-skeleton together with trivial weights all equal to $1$.
\end{example}

\begin{example}[A rigid curve] \label{example-rigidcurve}
    \begin{figure}[t]
      \includegraphics*[scale=0.4]{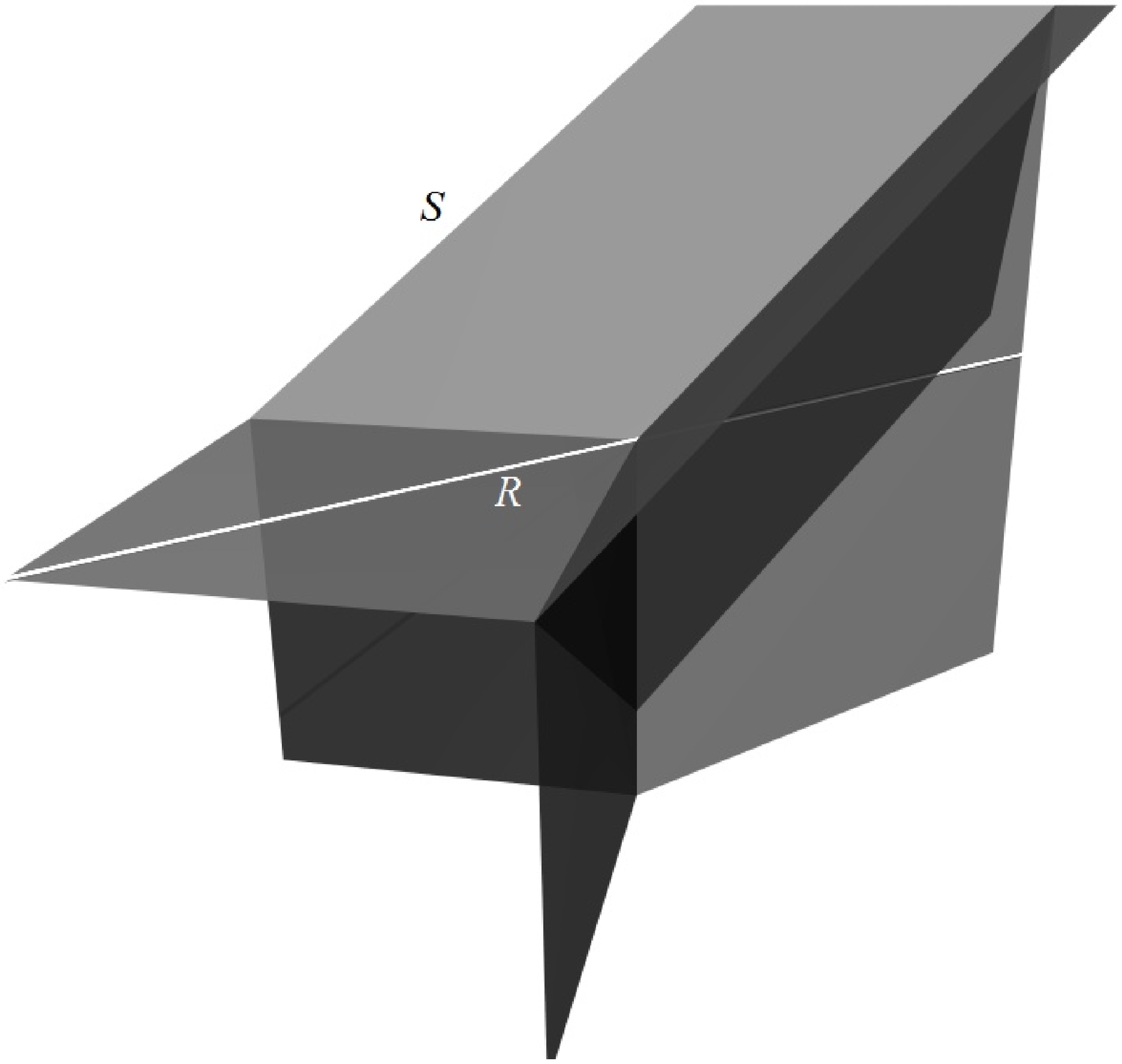}\\
      The rigid curve $R$ in $S$.
    \end{figure}
    Using notations from example \ref{example-selfintersectionofhyperplanes}, we consider as ambient cycle the surface $S := [L^3_2] = \calT(x_1 \oplus x_2 \oplus x_3 \oplus 0)$ in $\RR^3$. In this surface, we want to show that the curve $R := [(\RR \cdot e_R, \omega_R(\RR \cdot e_R) = 1)] \in Z_1^{\aff}(S)$, where $e_R := e_1 +e_2$, has negative self-intersection in the following sense: We construct a rational function $\varphi$ on $S$ whose associated Weil divisor is $R$ and show that $\varphi \cdot R = \varphi \cdot \varphi \cdot S$ is just the origin with weight $-1$. This curve and its rigidness were first discussed in \cite[Example 4.11., Example    5.9.]{M}.\\
    Let us construct $\varphi$. First we refine $L^3_2$ to $L_R$ by replacing $\sigma_{\{1,2\}}$ and $\sigma_{\{0,3\}}$ with $\sigma_{\{1,R\}}$, $\sigma_{\{R\}}$, $\sigma_{\{R,2\}}$, $\sigma_{\{0,-R\}}$, $\sigma_{\{-R\}}$ and $\sigma_{\{-R,3\}}$ (using again the notations from example \ref{example-selfintersectionofhyperplanes} and $e_{-R} := - e_R = e_0 + e_3$). We define $\varphi : |S| \rightarrow \RR$ to be the unique function that is linear on the faces of $L_R$ and fulfills
    $$
        0, -e_1, -e_2, -e_3, -e_{-R} \mapsto 0, \;\; -e_0 \mapsto 1 \text{ and } -e_R \mapsto -1.
    $$
    Analogous to \ref{example-selfintersectionofhyperplanes}, we can compute for $i = 1,2$
    $$
        \omega_{L_R, \varphi}(\sigma_{\{i\}}) = \varphi(-e_0) + \varphi(-e_3) + \varphi(-e_R) = 1 + 0 - 1 = 0,
    $$
    for $i = 0,3$
    $$
        \omega_{L_R, \varphi}(\sigma_{\{i\}}) = \varphi(-e_1) + \varphi(-e_2) + \varphi(-e_{-R}) = 0 + 0 + 0 = 0,
    $$
    and finally
    \begin{eqnarray*}
        & \omega_{L_R, \varphi}(\sigma_{\{R\}}) = \varphi(-e_1) + \varphi(-e_2) - \varphi(-e_R) = 0 + 0 -  (-1) = 1, & \\
        & \omega_{L_R, \varphi}(\sigma_{\{-R\}}) = \varphi(-e_0) + \varphi(-e_3) - \varphi(-e_{-R}) = 1 + 0 + 0 = 1, &
    \end{eqnarray*}
    which means $\varphi \cdot S = R$. Now we can easily compute $\varphi \cdot \varphi \cdot S = \varphi \cdot R$ on the representative $\{\sigma_{\{R\}}, \sigma_{\{-R\}}, \{0\} \}$ (with trivial weights) of $R$:
    $$
        \omega_{R, \varphi}(\{0\}) = \varphi(-e_R) + \varphi(-e_{-R}) = -1 + 0 = -1.
    $$
    Therefore $\varphi \cdot \varphi \cdot S = [(\{0\}, \omega(\{0\}) =
    -1)]$. Note that we really obtain a cycle with negative
    weight, not only a cycle class modulo rational equivalence as it is
    the case in ``classical'' algebraic geometry.
\end{example}

  \section{Push-forward of affine cycles and pull-back of Cartier divisors} \label{sec-pushpull}

The aim of this section is to construct push-forwards of cycles
and pull-backs of Cartier divisors along morphisms of fans and to
study the interaction of both constructions. To do this we first
of all have to introduce the notion of morphism:

\begin{definition}[Morphisms of fans] \label{morphism}
  Let $X$ be a fan in $V=\Lambda \otimes_\ZZ \RR$ and $Y$ be a fan
  in $V'=\Lambda' \otimes_\ZZ \RR$. A \df{morphism} $f: X \rightarrow Y$
  is a $\ZZ$-linear map, i.e. a map from $|X| \subseteq V$ to $|Y|
  \subseteq V'$ induced by a $\ZZ$-linear map $\widetilde{f}:\Lambda \rightarrow
  \Lambda'$. By abuse of notation we will usually denote all three maps $f,\widetilde{f}$ and
  $\widetilde{f} \otimes_{\ZZ} \text{id}$ by the same letter $f$ (note that
  the last two maps are in general not uniquely determined by $f:X \rightarrow Y$).
  A morphism of weighted fans is a morphism of fans. A morphism of
  affine cycles $f: [(X, \omega_X)] \rightarrow [(Y, \omega_Y)]$
  is a morphism of fans $f:X^{*} \rightarrow Y^{*}$.
  Note that in this latter case the notion of morphism does not depend on the choice
  of the representatives by remark \ref{rem-refinement}.
\end{definition}

Given such a morphism the following construction shows how to
build the push-forward fan of a given fan. Afterwards we will show
that this construction induces a well-defined operation on cycles.

\begin{construction} \label{constr-push}
  We refer to \cite[section 2]{GKM} for more details on the following construction.
  Let $(X,\omega_X)$ be a weighted fan of pure dimension $n$ in $V=\Lambda \otimes_\ZZ \RR$,
  let $Y$ be any fan in $V'=\Lambda' \otimes_\ZZ \RR$ and let $f:X \rightarrow Y$
  be a morphism. Passing to an appropriate refinement of $(X,\omega_X)$
  the collection of cones $$f_*X := \{ f(\sigma) | \sigma \in X \text{ contained
  in a maximal cone of } X \text{ on which } f \text{ is injective} \}$$
  is a fan in $V'$ of pure dimension $n$. It can be made into a
  weighted fan by setting $$\omega_{f_*X}(\sigma') := \sum_{\sigma \in X^{(n)}: f(\sigma)=\sigma'}
  \omega_X (\sigma) \cdot \left| \Lambda_{\sigma'}' / f(\Lambda_\sigma) \right|$$
  for all $\sigma' \in f_*X^{(n)}$. The equivalence class of this weighted fan only depends
  on the equivalence class of $(X,\omega_X)$.
\end{construction}

\begin{example}
  Let $X$ be the fan with cones $\tau_1, \tau_2, \tau_3, \{0\}$ as shown in the figure
  \medskip\\
  \begin{minipage}{\linewidth}
  \begin{center}
  \input{pic/ExPushf.pstex_t}\\
  \end{center}
  \end{minipage}
  \medskip\\
  and let $\omega_X(\tau_i)=1$ for $i=1,2,3$. Moreover, let $Y:=\RR$ be the real line and
  the morphisms $f_1, f_2: X \rightarrow Y$ be given by $f_1(x,y):=x+y$ and $f_2(x,y):=x$
  respectively. Then $(f_1)_{*}X = (f_2)_{*}X = \left\{ \{ x \leq 0\}, \{0\}, \{ x \geq 0\}  \right\}$, but $\omega_{(f_1)_{*}X}(\{ x \leq 0\})= {\omega_{(f_1)_{*}X}(\{ x \geq 0\})} =2$ and $\omega_{(f_2)_{*}X}(\{ x \leq 0\})= \omega_{(f_2)_{*}X}(\{ x \geq 0\}) =1$.
\end{example}

\begin{proposition} \label{prop-pushtropfans}
  Let $(X,\omega_X)$ be a tropical fan of dimension $n$ in $V=\Lambda \otimes_\ZZ \RR$,
  let $Y$ be any fan in $V'=\Lambda' \otimes_\ZZ \RR$ and let $f:X \rightarrow Y$
  be a morphism. Then $f_*X$ is a tropical fan of dimension $n$.
\end{proposition}
\begin{proof}
  A proof can be found in \cite[section 2]{GKM}.
\end{proof}

By construction \ref{constr-push} and proposition
\ref{prop-pushtropfans} the following definition is well-defined:

\begin{definition}[Push-forward of cycles]\label{def-pushcycles}
  Let $V=\Lambda \otimes_\ZZ \RR$ and $V'=\Lambda' \otimes_\ZZ
  \RR$. Moreover, let $X \in Z_m^{\aff}(V)$, $Y \in Z_n^{\aff}(V')$ and $f:X \rightarrow Y$ be a morphism.
  For $[(Z,\omega_Z)] \in Z_k^{\aff}(X)$ we define $$f_{*}[(Z,\omega_Z)] :=
  [(f_{*}(Z^{*}),\omega_{f_{*}(Z^{*})})] \in Z_k^{\aff}(Y).$$
\end{definition}

\begin{proposition}[Push-forward of cycles] \label{prop-pushcycles}
  Let $V=\Lambda \otimes_\ZZ \RR$ and $V'=\Lambda' \otimes_\ZZ
  \RR$. Let $X \in Z_m^{\aff}(V)$ and $Y \in Z_n^{\aff}(V')$ be cycles
  and let $f:X \rightarrow Y$ be a morphism. Then the map
  $$Z_k^{\aff}(X) \longrightarrow Z_k^{\aff}(Y): C \longmapsto f_{*}C$$
  is well-defined and $\ZZ$-linear.
\end{proposition}
\begin{proof}
  It remains to prove the linearity: Let $(A,\omega_A)$ and
  $(B,\omega_B)$ be two tropical fans of dimension $k$ with $A=A^{*}$, $B=B^{*}$ and
  $|A|,|B| \subseteq |X^{*}|$. We want to show that $f_{*}(A+B) \sim f_{*}A+f_{*}B$.
  Refining $A$ and $B$ as in construction \ref{constr-sum} we may assume that $A,B \subseteq
  A+B$. Set $\widetilde{A}:=A+B$ and
  $$\omega_{\widetilde{A}}(\sigma) := \left\{ \begin{array}{cl} \omega_{A} (\sigma), \text{ if }
  \sigma \in A \\ 0, \text{ else} \end{array} \right.$$
  for all facets $\sigma \in \widetilde{A}.$ Analogously, set $\widetilde{B} :=
  A+B$ with according weights.
  Then $\widetilde{A} \sim A$ and $\widetilde{B} \sim
  B$. Carrying out a further refinement of $A+B$ like in
  construction \ref{constr-push} we can reach that
  $f_{*}(A+B)=\{ f(\sigma) | {\sigma \in A+B} \text{ contained }
  \text{in a maximal cone of } A+B \text{ on which } f \text{ is injective}
  \}.$ Using $\widetilde{A}=\widetilde{B}=\widetilde{A}+\widetilde{B}=A+B$ we get
  $f_{*}\widetilde{A}=f_{*}\widetilde{B}={f_{*}(\widetilde{A}+\widetilde{B})}=f_{*}(A+B)$
  and it remains to compare the weights:
  {\begin{eqnarray*}
    \omega_{f_{*}(\widetilde{A}+\widetilde{B})}(\sigma')&=& \sum_{\sigma \in (\widetilde{A}+\widetilde{B})^{(k)}: f(\sigma)=\sigma'}
    \omega_{\widetilde{A}+\widetilde{B}} (\sigma) \cdot \left| \Lambda_{\sigma'}' / f(\Lambda_\sigma) \right| \\
    &=& \sum_{\sigma \in (\widetilde{A}+\widetilde{B})^{(k)}: f(\sigma)=\sigma'}
    \left[ \omega_{\widetilde{A}} (\sigma) + \omega_{\widetilde{B}} (\sigma) \right] \cdot \left| \Lambda_{\sigma'}' / f(\Lambda_\sigma) \right|\\
    &=& \sum_{\sigma \in \widetilde{A}^{(k)}: f(\sigma)=\sigma'}
    \omega_{\widetilde{A}} (\sigma) \cdot \left| \Lambda_{\sigma'}' / f(\Lambda_\sigma)
    \right| + \\
    && \sum_{\sigma \in \widetilde{B}^{(k)}: f(\sigma)=\sigma'}
    \omega_{\widetilde{B}} (\sigma) \cdot \left| \Lambda_{\sigma'}' / f(\Lambda_\sigma) \right|\\
    &=& \omega_{f_{*}\widetilde{A}}(\sigma')+\omega_{f_{*}\widetilde{B}}(\sigma')
  \end{eqnarray*}}
  for all facets $\sigma'$ of $f_{*}(A+B)$.
  Hence $f_{*}(A+B) \sim f_{*}(\widetilde{A}+\widetilde{B}) =
  f_{*}\widetilde{A}+f_{*}\widetilde{B} \sim f_{*}A+f_{*}B$ as weighted fans.
\end{proof}

Our next step is now to define the pull-back of a Cartier divisor.
As promised we will prove after this a projection formula that
describes the interaction between our two constructions.

\begin{proposition}[Pull-back of Cartier divisors] \label{prop-pulldivisors}
  Let $C \in Z_m^{\aff}(V)$ and $D \in Z_n^{\aff}(V')$ be cycles in $V=\Lambda \otimes_\ZZ \RR$ and
  $V'=\Lambda' \otimes_\ZZ \RR$ respectively and let $f:C \rightarrow D$
  be a morphism. Then there is a well-defined and $\ZZ$-linear map
  $$\Div(D) \longrightarrow \Div(C): [h] \longmapsto f^{*}[h]:=[h \circ f].$$
\end{proposition}
\begin{proof}
  The map $h \mapsto h \circ f$ is obviously $\ZZ$-linear on rational functions and maps
  affine linear functions to affine linear functions. Thus it remains to prove that $h \circ f$
  is a rational function if $h$ is one:
  Therefore let $(X,\omega_X)$ be any representative of $C$, let $(Y,\omega_Y)$ be a reduced representative
  of $D$ such that the restriction of $h$ to every cone in $Y$ is affine linear and let $f_V: V \rightarrow V'$
  be a $\ZZ$-linear map such that $f_V |_{|C|}=f$.
  Since $Z:= \linebreak \{f_V^{-1}(\sigma')|\sigma' \in Y\}$ is a fan in $V$ and
  $|X| \subseteq |Z|$ we can construct the refinement $\widetilde{X} := X \cap Z$ of
  $X$ such that $h \circ f$ is affine linear on every cone of
  $\widetilde{X}$. This finishes the proof.
\end{proof}

\begin{proposition}[Projection formula] \label{prop-projectionf}
  Let $C \in Z_m^{\aff}(V)$ and $D \in Z_n^{\aff}(V')$ be cycles in $V=\Lambda \otimes_\ZZ \RR$ and
  $V'=\Lambda' \otimes_\ZZ \RR$ respectively and let $f:C \rightarrow D$
  be a morphism. Let $E \in Z_k^{\aff}(C)$ be a cycle and
  let $\varphi \in \Div(D)$ be a Cartier divisor. Then the following equation holds:
  $$ \varphi \cdot (f_{*}E) = f_{*}(f^{*}\varphi \cdot E) \in Z_{k-1}^{\aff}(D).$$
\end{proposition}
\begin{proof}
  Let $E=[(Z,\omega_Z)]$ and $\varphi=[h]$. We may assume that $Z=Z^{*}$ and $h(0)=0$. Replacing $Z$ by a refinement we may
  additionally assume that $f^{*}h$ is linear on every cone
  of $Z$ (cf. definition \ref{def-rationalfunctions}) and that $$f_{*}Z = \{f(\sigma)| \sigma \in Z \text{ contained
  in a maximal cone of } Z \text{ on which } f \text{ is injective}
  \}$$ (cf. construction \ref{constr-push}). Note that in this
  case $h$ is linear on the cones of $f_{*}Z$, too. Let ${\sigma' \subseteq |D|}$ be a cone
  (not necessarily $\sigma' \in f_{*}Z$) such that $h$ is linear on $\sigma'$. Then there is
  a unique linear map $h_{\sigma'}: V_{\sigma'}' \rightarrow \RR$
  induced by the restriction $h|_{\sigma'}$. Analogously for
  $f^{*}h_{\sigma}, \sigma \subseteq |C|$.
  For cones $\tau<\sigma \in Z$ of dimension $k-1$ and $k$ respectively
  let $v_{\sigma / \tau} \in \Lambda$ be a representative of the
  primitive normal vector $u_{\sigma / \tau} \in \Lambda / \Lambda_{\tau}$ of
  construction \ref{constr-normalvectors}.
  Analogously, for $\tau'< \sigma' \in f_{*}Z$ of dimension $k-1$ and $k$ respectively let
  $v_{\sigma'/\tau'}$ be a representative of $u_{\sigma'/\tau'} \in \Lambda'/\Lambda_{\tau'}'$.
  Now we want to compare the weighted fans $h \cdot (f_{*}Z)$ and $f_{*}(f^{*}h \cdot Z)$:
  Let $\tau' \in f_{*}Z$ be a cone of dimension $k-1$.
  Then we can calculate the weight of $\tau'$ in $h \cdot (f_{*}Z)$ as follows:
  {\small \begin{eqnarray*}
    \omega_{h \cdot (f_{*}Z)} (\tau') &=& \left( \sum_{\sigma' \in f_{*}Z: \sigma' >
    \tau'} \omega_{f_{*}Z}(\sigma') \cdot h_{\sigma'}(v_{\sigma'/\tau'})
    \right)\\
    && -h_{\tau'}
    \left( \sum_{\sigma' \in f_{*}Z: \sigma' > \tau'} \omega_{f_{*}Z}(\sigma')
    \cdot v_{\sigma'/\tau'} \right)\\
    &=& \left( \sum_{\sigma' \in f_{*}Z: \sigma' > \tau'}
    \left( \sum_{\sigma \in Z^{(k)}:f(\sigma)=\sigma'} \omega_Z(\sigma) \cdot | \Lambda_{\sigma'}'/f(\Lambda_{\sigma}) | \right)
    \cdot h_{\sigma'}(v_{\sigma'/\tau'}) \right)\\
    && - h_{\tau'} \left( \sum_{\sigma' \in f_{*}Z: \sigma' > \tau'} \left( \sum_{\sigma \in Z^{(k)}:f(\sigma)=\sigma'} \omega_Z(\sigma) \cdot | \Lambda_{\sigma'}'/f(\Lambda_{\sigma}) | \right)
    \cdot v_{\sigma'/\tau'} \right)\\
    &=& \left( \sum_{\sigma \in Z^{(k)}:f(\sigma)>\tau'} \omega_Z(\sigma) \cdot | \Lambda_{f(\sigma)}'/f(\Lambda_{\sigma})
    | \cdot h_{f(\sigma)}(v_{f(\sigma)/\tau'}) \right)\\ &&
    - h_{\tau'}
    \left( \sum_{\sigma \in Z^{(k)}:f(\sigma)>\tau'} \omega_Z(\sigma) \cdot | \Lambda_{f(\sigma)}'/f(\Lambda_{\sigma})
    | \cdot v_{f(\sigma)/\tau'} \right)
  \end{eqnarray*}}
  Now let $\tau' \in f_{*}(f^{*}h \cdot Z)$ of dimension $k-1$.
  The weight of $\tau'$ in $f_{*}(f^{*}h \cdot Z)$ can be calculated as follows:
  {\small \begin{eqnarray*}
    \omega_{f_{*}(f^{*}h \cdot Z)}(\tau') &=& \sum_{\twolines{\tau \in (f^{*}h \cdot Z)^{(k-1)}:}{f(\tau)=\tau'}}
    \omega_{f^{*}h \cdot Z}(\tau) \cdot | \Lambda_{\tau'}' / f(\Lambda_{\tau})|\\
    &=& \sum_{\twolines{\tau \in (f^{*}h \cdot
    Z)^{(k-1)}:}{f(\tau)=\tau'}}
    \left( \sum_{\sigma \in Z^{(k)}: \sigma>\tau} \omega_Z(\sigma) f^{*}h_{\sigma}(v_{\sigma/\tau}) \right.\\
    && - \left. f^{*}h_{\tau} \left( \sum_{\sigma \in Z^{(k)}: \sigma>\tau} \omega_Z(\sigma) \cdot v_{\sigma/\tau} \right) \right) \cdot | \Lambda_{\tau'}' / f(\Lambda_{\tau})|
  \end{eqnarray*}}
  {\small \begin{eqnarray*}
    \phantom{\omega_{f_{*}(f^{*}h \cdot Z)}(\tau')}
    &=& \sum_{\twolines{\tau \in (f^{*}h \cdot Z)^{(k-1)}:}{f(\tau)=\tau'}}
    \left( \sum_{\sigma \in Z^{(k)}: \sigma>\tau} \omega_Z(\sigma) h_{f(\sigma)}(f(v_{\sigma/\tau})) \right. \\
    && - \left. h_{f(\tau)} \left( \sum_{\sigma \in Z^{(k)}: \sigma>\tau} \omega_Z(\sigma) \cdot f \left( v_{\sigma/\tau} \right) \right) \right) \cdot | \Lambda_{\tau'}' / f(\Lambda_{\tau})
    |.\\
  \end{eqnarray*}}
  Note that $f(v_{\sigma/\tau})= {| \Lambda_{\sigma'}'/(\Lambda_{\tau'}'+
  \ZZ f(v_{\sigma/\tau}))|} \cdot v_{\sigma' / \tau'} +
  \lambda_{\sigma,\tau} \in \Lambda'$ for some $\lambda_{\sigma,\tau} \in
  \Lambda_{\tau'}'$. Since $h_{f(\sigma)}(\lambda_{\sigma,\tau})=h_{f(\tau)}(\lambda_{\sigma,\tau})$
  these parts of the corresponding summands in the first and second interior sum cancel using the linearity of $h_{f(\tau)}$.
  Moreover, note that $f(v_{\sigma/\tau}) =\lambda_{\sigma,\tau} \in \Lambda_{\tau'}'$ for those $\sigma>\tau$ on which $f$
  is not injective and that the whole summands cancel in this case. Thus we can conclude that the sum does not
  change if we restrict the summation to those $\sigma > \tau$ on which $f$ is injective.
  Using additionally the equation $${|\Lambda_{\sigma'}'/f(\Lambda_{\sigma}) |} =
  {|\Lambda_{\tau'}'/f(\Lambda_{\tau})| \cdot
  |\Lambda_{\sigma'}'/(\Lambda_{\tau'}'+ \ZZ
  f(v_{\sigma/\tau}))|}$$ we get
  {\small \begin{eqnarray*}
    \omega_{f_{*}(f^{*}h \cdot Z)}(\tau')
    &=& \sum_{\twolines{\tau \in (f^{*}h \cdot Z)^{(k-1)}:}{f(\tau)=\tau'}}
    \left( \sum_{\twolines{\sigma \in Z^{(k)}:}{\sigma>\tau,f(\sigma)>\tau'}} \omega_Z(\sigma) \cdot | \Lambda_{f(\sigma)}'/f(\Lambda_{\sigma})| \cdot h_{f(\sigma)}(v_{f(\sigma)/\tau'}) \right. \\
    && - \left. h_{\tau'} \left( \sum_{\twolines{\sigma \in Z^{(k)}:}{\sigma>\tau,f(\sigma)>\tau'}} \omega_Z(\sigma) \cdot | \Lambda_{f(\sigma)}'/f(\Lambda_{\sigma})| \cdot v_{f(\sigma)/\tau'}
    \right) \right)\\
    &=& \left( \sum_{\sigma \in Z^{(k)}: f(\sigma)>\tau'}
    \omega_Z(\sigma) \cdot | \Lambda_{f(\sigma)}'/f(\Lambda_{\sigma})| \cdot h_{f(\sigma)}(v_{f(\sigma)/\tau'}) \right) \\
    && - h_{\tau'} \left( \sum_{\sigma \in Z^{(k)}: f(\sigma)>\tau'} \omega_Z(\sigma) \cdot | \Lambda_{f(\sigma)}'/f(\Lambda_{\sigma})| \cdot
    v_{f(\sigma)/\tau'} \right).
  \end{eqnarray*}}
  Note that for the last equation we used again the linearity
  of $h_{\tau'}$. We have checked so far that a cone $\tau'$ of dimension $k-1$
  occurring in both $h \cdot (f_{*}Z)$ and $f_{*}(f^{*}h \cdot Z)$ has the
  same weight in both fans. Thus it remains to examine those cones
  $f(\tau), \tau \in Z^{(k-1)}$ such that $f$ is injective on $\tau$ but not on any
  $\sigma > \tau$: In this case all vectors $v_{\sigma/\tau}$ are
  mapped to $\Lambda_{f(\tau)}'$. Again, $h_{f(\sigma)}=h_{f(\tau)}$ and by linearity of
  $h_{f(\tau)}$ all summands in the sum cancel as above.
  Hence the the weight of $f(\tau)$ in $f_{*}(f^{*}h \cdot Z)$ is 0 and
  $\varphi \cdot (f_{*}E)=[h \cdot (f_{*}Z)] = [f_{*}(f^{*}h \cdot Z)]=f_{*}(f^{*} \varphi \cdot E)$.
\end{proof}

  \section{Abstract tropical cycles} \label{sec-cycles}

In this section we will introduce the notion of abstract tropical
cycles as spaces that have tropical fans as local building blocks.
Then we will generalize the theory from the previous sections to
these spaces.

\begin{definition}[Abstract polyhedral complexes] \label{defn-polycomplex}
  An \df{(abstract) polyhedral complex} is a topological space
  $|X|$ together with a finite set $X$ of closed subsets of $|X|$
  and an embedding map $\varphi_{\sigma}: \sigma \rightarrow
  \RR^{n_{\sigma}}$ for every $\sigma \in X$ such that
  \begin{enumerate}
    \item $X$ is closed under taking intersections, i.e. $\sigma \cap \sigma' \in X$
    for all $\sigma, \sigma' \in X$ with $\sigma \cap \sigma' \neq
    \emptyset$,
    \item every image $\varphi_{\sigma}(\sigma)$, $\sigma \in X$ is a rational
    polyhedron not contained in a proper affine subspace of $\RR^{n_{\sigma}}$,
    \item for every pair $\sigma, \sigma' \in X$
    the concatenation $\varphi_{\sigma} \circ \varphi_{\sigma'}^{-1}$ is integer
    affine linear where defined,
    \item $|X|=\bigcup\limits_{\sigma \in X}^{\mbox{\Large .}}
    \varphi_{\sigma}^{-1}(\varphi_{\sigma}(\sigma)°)$, where
    $\varphi_{\sigma}(\sigma)°$ denotes the interior of
    $\varphi_{\sigma}(\sigma)$ in $\RR^{n_{\sigma}}$.
  \end{enumerate}
  For simplicity we will usually drop the embedding maps $\varphi_{\sigma}$
  and denote the polyhedral complex $(X,|X|,\{ \varphi_{\sigma} | \sigma \in X
  \})$ by $(X,|X|)$ or just by $X$ if no confusion can occur.
  The closed subsets $\sigma \in X$ are called the \df{polyhedra}
  or \df{faces of $(X,|X|)$}. For $\sigma \in X$ the open set $\sigma^{ri} :=
  \varphi_{\sigma}^{-1}(\varphi_{\sigma}(\sigma)°)$ is called the
  \df{relative interior of $\sigma$}. Like in the case of fans the \df{dimension} of
  $(X,|X|)$ is the maximum of the dimensions of its polyhedra.
  $(X,|X|)$ is \df{pure-dimensional} if every inclusion-maximal
  polyhedron has the same dimension. We denote by $X^{(n)}$ the set
  of polyhedra in $(X,|X|)$ of dimension $n$. Let $\tau, \sigma \in X$. Like in the case of
  fans we write $\tau \leq \sigma$ (or $\tau < \sigma$) if $\tau \subseteq \sigma$
  (or $\tau \subsetneq \sigma$ respectively). \medskip \\
  An abstract polyhedral complex $(X,|X|)$ of pure dimension $n$
  together with a map $\omega_X : X^{(n)} \rightarrow \ZZ$ is called
  \df{weighted polyhedral complex} of dimension $n$ and
  $\omega_X(\sigma)$ the \df{weight} of the polyhedron $\sigma \in X^{(n)}$.
  Like in the case of fans the empty complex $\emptyset$ is a
  weighted polyhedral complex of every dimension $n$.
  If $\left( (X,|X|), \omega_X \right)$ is a weighted polyhedral
  complex of dimension $n$ then let $$X^{*}:= {\{ \tau \in X | \tau \subseteq \sigma \text{ for
  some } \sigma \in X^{(n)} \text { with } \omega_X(\sigma) \neq 0
  \}}, |X^{*}| := \bigcup_{\tau \in X^{*}} \tau \subseteq |X|.$$
  With these definitions $\left( (X^{*},|X^{*}|), \omega_X|_{(X^{*})^{(n)}} \right)$
  is again a weighted polyhedral complex of dimension $n$,
  called the \df{non-zero part} of $\left( (X,|X|), \omega_X \right)$.
  We call a weighted polyhedral complex $\left( (X,|X|), \omega_X \right)$
  \df{reduced} if $\left( (X,|X|), \omega_X \right)=\left( (X^{*},|X^{*}|), \omega_{X^{*}} \right)$
  holds.
\end{definition}

\begin{definition}[Subcomplexes and refinements] \label{defn-subcomplex}
  Let $(X,|X|,\{ \varphi_{\sigma}\})$ and $(Y,|Y|,\{ \psi_{\tau}\})$
  be two polyhedral complexes. We call $(X,|X|,\{ \varphi_{\sigma}\})$
  a \df{subcomplex} of $(Y,|Y|,\{ \psi_{\tau}\})$ if
  \begin{enumerate}
    \item $|X| \subseteq |Y|$,
    \item \label{pb} for every $\sigma \in X$ exists $\tau \in Y$ with $\sigma \subseteq \tau$ and
    \item \label{pc} the $\ZZ$-linear structures of $X$ and $Y$ are compatible, i.e. for a
    pair $\sigma, \tau$ from \ref{pb} the maps $\varphi_{\sigma} \circ \psi_{\tau}^{-1}$ and $\psi_{\tau} \circ \varphi_{\sigma}^{-1}$
    are integer affine linear where defined.
  \end{enumerate}
  We write $(X,|X|,\{ \varphi_{\sigma}\}) \subfan (Y,|Y|,\{ \psi_{\tau}\})$ in this case.
  Analogous to the case of fans we define a map $C_{X,Y} : X
  \rightarrow Y$ that maps a polyhedron in $X$ to the
  inclusion-minimal polyhedron in $Y$ containing it.\\
  We call a weighted polyhedral complex $\left( (X,|X|),\omega_X \right)$
  a \df{refinement} of $\left( (Y,|Y|),\omega_Y \right)$ if
  \begin{enumerate}
    \item $(X^{*},|X^{*}|) \subfan (Y^{*},|Y^{*}|)$,
    \item $|X^{*}| = |Y^{*}|$,
    \item $\omega_X(\sigma) = \omega_Y(C_{X^{*},Y^{*}}(\sigma))$ for all
    $\sigma \in (X^{*})^{(\dim(X))}$.
  \end{enumerate}
\end{definition}

\begin{definition}[Open fans] \label{defn-openfan}
  Let $(\widetilde{F},\omega_{\widetilde{F}})$ be a tropical fan in
  $\RR^n$ and $U \subseteq \RR^n$ an open subset containing the origin. The set
  $F:= \widetilde{F} \cap U := \{\sigma \cap U | \sigma \in \widetilde{F} \}$ together with the induced
  weight function $\omega_F$ is called an \df{open (tropical) fan} in
  $\RR^n$. Like in the case of fans let $|F|:= \bigcup_{\sigma' \in F}
  \sigma'$. Note that the open fan $F$ contains the whole
  information of the entire fan $\widetilde{F}$ as
  $\widetilde{F}= \{ \RR_{\geq 0} \cdot \sigma' | \sigma' \in F
  \}$.
\end{definition}

\begin{definition}[Tropical polyhedral complexes] \label{defn-tropcomplex}
  A \df{tropical polyhedral complex} of dimension $n$ is a weighted polyhedral complex
  $\left( (X,|X|), \omega_X \right)$ of pure dimension $n$
  together with the following data:
  For every polyhedron $\sigma \in X^{*}$ we are given an open fan
  $F_\sigma$ in some $\RR^{n_\sigma}$ and a homeomorphism
  $$\Phi_{\sigma}: S_{\sigma} := \bigcup_{\sigma' \in X^{*}, \sigma' \supseteq \sigma}
  (\sigma')^{ri} \stackrel{\sim}{\longrightarrow} |F_{\sigma}|$$
  such that
  \begin{enumerate}
    \item for all $\sigma' \in X^{*}, \sigma' \supseteq \sigma$ holds $\Phi_{\sigma}(\sigma' \cap S_{\sigma}) \in F_{\sigma}$
    and $\Phi_{\sigma}$ is compatible with the $\ZZ$-linear structure
    on $\sigma'$, i.e. $\Phi_{\sigma} \circ \varphi_{\sigma'}^{-1}$
    and $\varphi_{\sigma'} \circ \Phi_{\sigma}^{-1}$ are integer affine
    linear where defined,
    \item $\omega_X(\sigma') = \omega_{F_\sigma} (\Phi_{\sigma}
    (\sigma' \cap S_{\sigma}))$ for every $\sigma' \in (X^{*})^{(n)}$ with
    $\sigma' \supseteq \sigma$,
    \item for every pair $\sigma, \tau \in X^{*}$ there is an integer affine linear map $A_{\sigma,\tau}$
    and a commutative diagram
    $$\begin{xy}
      \xymatrix{ S{_\sigma} \cap S_{\tau} \ar[d]^{\sim}_{\Phi_{\sigma}} \ar[r]_{\sim}^{\Phi_{\tau}} &
      \Phi_{\tau}(S_{\sigma} \cap S_{\tau})\\
      \Phi_{\sigma}(S_{\sigma} \cap S_{\tau}) \ar[ur]_{A_{\sigma,\tau}} }
    \end{xy}.$$
  \end{enumerate}
  For simplicity of notation we will usually drop the maps $\Phi_{\sigma}$
  and write $\left( (X,|X|), \omega_X \right)$ or just $X$ instead of $\left( ((X,|X|), \omega_X),
  \{\Phi_{\sigma} \} \right)$.
  A tropical polyhedral complex is called \df{reduced} if the
  underlying weighted polyhedral complex is.
\end{definition}

\begin{example} \label{example-abstractcycles}
  The following figure shows the topological spaces and the decompositions into polyhedra of two such abstract tropical polyhedral complexes
  together with the open fan $F_\sigma$ for every polyhedron $\sigma$:
  \medskip\\
  \begin{minipage}{\linewidth}
  \begin{center}
  \input{pic/ExTropComplex.pstex_t}\\
  \end{center}
  \end{minipage}
\end{example}

\begin{construction}[Refinements of tropical polyhedral complexes] \label{constr-refinementtropcomplex}
  Let \linebreak $\left( ((X,|X|), \omega_X), \{\Phi_\sigma\} \right)$ be a tropical polyhedral complex
  and let $((Y,|Y|), \omega_Y)$ be a refinement of its underlying weighted polyhedral complex $((X,|X|), \omega_X)$.
  Then we can make $((Y,|Y|), \omega_Y)$ into a tropical polyhedral complex as follows:
  We may assume that $X$ and $Y$ are reduced
  as we do not pose any conditions on polyhedra with weight zero. Fix
  some  $\tau \in Y$ and let $\sigma:= C_{Y,X}(\tau)$. By definition of refinement,
  for every $\tau' \in Y$ with $\tau' \geq \tau$ there is
  $\sigma' \in X$, $\sigma' \geq \sigma$ with $\tau' \subseteq \sigma'$.
  Thus $S_{\tau} \subseteq S_{\sigma}$ and we have a map $\Psi_{\tau}:=\Phi_{\sigma}|_{S_{\tau}} : S_{\tau}
  \stackrel{\sim}{\rightarrow} \Psi_{\tau}(S_{\tau}) \subseteq
  \RR^{n_{\sigma}}$. It remains to give $\Psi_{\tau}(S_{\tau})$
  the structure of an open fan: We may assume that
  $\{0\} \subseteq \Psi_{\tau}(\tau)$ (otherwise replace $\Psi_{\tau}$ by the concatenating
  of $\Psi_{\tau}$ with an appropriate translation $T_{\tau}$, apply $T_{\tau}$ to $F_{\sigma}^{X}$ and $\Phi_{\sigma}$
  and change the maps $A_{\sigma,\sigma'}$ and $A_{\sigma',\sigma}$ accordingly).
  Let $\widetilde{F}_{\sigma}^{X}:= {\{ \RR_{\geq 0} \cdot \sigma' | \sigma' \in
  F_{\sigma}^{X} \}}$ be the tropical fan associated to
  $F_{\sigma}^{X}$ and let $\widetilde{F}_{\tau}^{Y}$ be the set of
  cones $\widetilde{F}_{\tau}^{Y} := {\{ \RR_{\geq 0} \cdot \Psi_{\tau}(\tau') | \tau \leq
  \tau' \in Y \}}$. Note that the conditions on the $\ZZ$-linear structures on $X$ and $Y$ to be compatible
  and on $\Phi_{\sigma}$ to be compatible with the $\ZZ$-linear structure on $X$
  assure that $\widetilde{F}_{\tau}^{Y}$ is a fan in $\RR^{n_{\sigma}}$.
  In fact, $\widetilde{F}_{\tau}^{Y}$ with the weights induced by $Y$ is a refinement of
  $(\widetilde{F}_{\sigma}^{X},\omega_{\widetilde{F}_{\sigma}^{X}})$. Thus the maps $\Psi_{\tau}$ together
  with the open fans $\{ \varrho \cap \Psi_{\tau}(S_{\tau})| \varrho \in
  \widetilde{F}_{\tau}^{Y}\}$, $\tau \in Y$ fulfill all requirements for a tropical polyhedral
  complex.
\end{construction}

\begin{remark} \label{rem-constrtroprefinem}
  If not stated otherwise we will from now on equip every
  refinement of a tropical polyhedral complex coming from a
  refinement of the underlying weighted polyhedral complex with
  the tropical structure constructed in
  \ref{constr-refinementtropcomplex}.
\end{remark}

\begin{definition}[Refinements and equivalence of tropical polyhedral complexes] \label{def-equivrelontropcycles}
  Let $C_1=\left( ((X_1,|X_1|), \omega_{X_1}), \{\Phi_{\sigma_1}^{X_1}\} \right)$ and $C_2=\left( ((X_2,|X_2|),
  \omega_{X_2}), \{\Phi_{\sigma_2}^{X_2}\} \right)$ be tropical polyhedral complexes. We call $C_2$ a \df{refinement}
  of $C_1$ if
  \begin{enumerate}
    \item $((X_2,|X_2|),\omega_{X_2})$ is a refinement of $((X_1,|X_1|), \omega_{X_1})$ and
    \item \label{lb} $C_2$ carries the tropical structure induced by $C_1$ like in construction
    \ref{constr-refinementtropcomplex}, i.e. if $C_2'= \left( ((X_2,|X_2|), \omega_{X_2}), \{\widetilde{\Phi}_{\sigma_2}^{X_2}\}
    \right)$ is the tropical polyhedral complex obtained from $C_1$ and the refinement $((X_2,|X_2|), \omega_{X_2})$
    then the maps ${\widetilde{\Phi}_{\sigma_2}^{X_2} \circ (\Phi_{\sigma_2}^{X_2})^{-1}}$ and ${\Phi_{\sigma_2}^{X_2} \circ
    (\widetilde{\Phi}_{\sigma_2}^{X_2})^{-1}}$ are integer affine linear where defined.
  \end{enumerate}
  We call two tropical polyhedral complexes $C_1$ and $C_2$ \df{equivalent} (write $C_1 \sim C_2$) if
  they have a common refinement (as tropical polyhedral complexes).
\end{definition}

\begin{remark} \label{rem-reftroppolycomplex}
  Note that different choices of translation maps $T_{\tau}$ in construction \ref{constr-refinementtropcomplex}
  only lead to tropical polyhedral complexes carrying the same tropical structure in the sense of
  definition \ref{def-equivrelontropcycles} \ref{lb}. In particular definition \ref{def-equivrelontropcycles} does
  not depend on the choices we made in construction \ref{constr-refinementtropcomplex}.
  Note moreover that refinements of $\left( ((X,|X|), \omega_X), \{\Phi_\sigma\} \right)$ and $((Y,|Y|), \omega_Y)$
  in construction \ref{constr-refinementtropcomplex} only lead to refinements of $\left( ((Y,|Y|),\omega_Y), \{\Psi_{\tau}\} \right)$.
\end{remark}

\begin{construction}(Refinements) \label{constr-constrrefinemofwpolycompl}
  Let $\left(((X,|X|,\{\varphi_{\sigma}\}),\omega_X),\{ \Phi_{\sigma} \}\right)$
  and \linebreak $\left(((Y,|Y|,\{\psi_{\tau}\}),\omega_Y),\{ \Psi_{\tau}\}\right)$ be
  reduced tropical polyhedral complexes such that \linebreak $(Y,|Y|) \subfan (X,|X|)$ and
  the tropical structures on $X$ and $Y$ agree, i.e. for every $\tau \in Y$ and
  $\sigma:= C_{Y,X}(\tau) \in X$ the maps $\Psi_{\tau} \circ \Phi_{\sigma}^{-1}$ and
  $\Phi_{\sigma} \circ \Psi_{\tau}^{-1}$ are integer affine linear where defined. Moreover let
  $\left(((X',|X'|,\{\varphi_{\sigma'}'\}), \omega_{X'}),\{ \Phi_{\sigma'}' \}\right)$
  be a reduced \linebreak refinement of $\left(((X,|X|,\{\varphi_{\sigma}\}),\omega_X),\{ \Phi_{\sigma} \}\right)$.
  Like in the case of fans we will \linebreak construct a refinement $\left((({Y \cap X'}, {|Y \cap X'|},
  \{\psi_{\tau'}^{Y \cap X'}\}), \omega_{Y \cap X'}), \{ \Psi_{\tau'}^{Y \cap X'}\} \right)$
  of \linebreak $\left(((Y,|Y|,\{\psi_{\tau}\}),\omega_Y),\{ \Psi_{\tau}\}\right)$
  such that $({Y \cap X'}, {|Y \cap X'|}) \subfan (X',|X'|)$ and the tropical structures
  on $Y \cap X'$ and $X'$ agree:\\
  Fix $\sigma \in X$. Note that the compatibility conditions on
  the $\ZZ$-linear structures of $X'$, $X$ and $Y$, $X$
  respectively (cf. \ref{defn-subcomplex} \ref{pc}) assure that
  $\varphi_{\sigma}(\sigma')$, $\sigma' \in X'$ with $\sigma' \subseteq \sigma$ as well as
  $\varphi_{\sigma}(\tau)$, $\tau \in Y$ with $\tau \subseteq \sigma$
  are rational polyhedra in $\RR^{n_{\sigma}}$. Thus in this case
  $\varphi_{\sigma}(\sigma' \cap \tau)=\varphi_{\sigma}(\sigma') \cap \varphi_{\sigma}(\tau)$
  is a rational polyhedron, too. Let $H_{\sigma',\tau} \cong \RR^{n_{\tau}}$ be the smallest affine subspace
  of $\RR^{n_{\sigma}}$ containing $\varphi_{\sigma}(\sigma' \cap \tau)$. We can consider
  $\varphi_{\sigma}|_{\sigma' \cap \tau}$ to be a map $\sigma' \cap \tau \rightarrow
  \RR^{n_{\tau}}$. We can hence construct the underlying weighted
  polyhedral complex of our desired tropical polyhedral complex
  as follows: Set $P:= \{ \tau \cap \sigma' | \tau \in Y, \sigma' \in X'\}$,
  $Y \cap X':= \{ \tau \in P| \nexists \widetilde{\tau} \in P^{(\dim(\tau))}: \widetilde{\tau} \subsetneq \tau \}$,
  $|Y \cap X'|:=|Y|$ and $\omega_{Y \cap X'}(\tau) := \omega_{Y}(C_{Y \cap X',Y}(\tau))$
  for all $\tau \in (Y \cap X')^{(\dim(Y))}$. It remains to define the
  maps $\psi_{\tau'}^{Y \cap X'}$ and $\Psi_{\tau'}^{Y \cap X'}$: For every $\tau' \in Y \cap X'$ choose a triplet
  $\sigma' \in X', \tau \in Y, \sigma \in X$ such that $\sigma'
  \cap \tau = \tau'$ and $\sigma', \tau \subseteq \sigma$ and set
  $\psi_{\tau'}^{Y \cap X'}:=\varphi_{\sigma}|_{\sigma' \cap \tau}$. With
  these definitions the weighted polyhedral complex $(({Y \cap X'}, {|Y \cap X'|},
  \{\psi_{\tau'}^{Y \cap X'}\}), \omega_{Y \cap X'})$
  is a refinement of $((Y,|Y|,\{\psi_{\tau}\}),\omega_Y)$. Thus we
  can apply construction \ref{constr-refinementtropcomplex} to
  obtain maps $\{ \Psi_{\tau'}^{Y \cap X'} \}$ that endow our
  weighted polyhedral complex with the tropical structure
  inherited from $((Y,|Y|,\{\psi_{\tau}\}),\omega_Y)$. Note that
  the compatibility property between the tropical structures of
  $Y$ and $X$ is bequeathed to $Y \cap X'$ and $X'$, too.
\end{construction}

\begin{lemma} \label{lemma-equivrelontroppolycomplex}
  The equivalence of tropical polyhedral complexes is an equivalence relation.
\end{lemma}
\begin{proof}
  Let $C_1=\left( ((X_1,|X_1|), \omega_{X_1}), \{\Phi_{\sigma_1}^{X_1}\} \right)$, $C_2=\left( ((X_2,|X_2|),
  \omega_{X_2}), \{\Phi_{\sigma_2}^{X_2}\} \right)$ and $C_3=\left( ((X_3,|X_3|), \omega_{X_3}),
  \{\Phi_{\sigma_3}^{X_3}\} \right)$ be tropical polyhedral complexes such that $C_1 \sim C_2$ via
  a common refinement $D_1=\left( ((Y_1,|Y_1|), \omega_{Y_1}), \{\Phi_{\sigma_1}^{Y_1}\}
  \right)$ and $C_2 \sim C_3$ via a common refinement $D_2=\left( ((Y_2,|Y_2|), \omega_{Y_2}), \{\Phi_{\sigma_2}^{Y_2}\}
  \right)$. We have to construct a common refinement of $C_1$ and $C_3$:
  First of all we may assume that $D_1$ and $D_2$ are reduced.
  Using construction \ref{constr-constrrefinemofwpolycompl} we get a refinement
  $D_3 := \left( ((Y_1 \cap Y_2,|Y_1 \cap Y_2|),\omega_{Y_1 \cap Y_2}), \{ \Phi_{\tau}^{Y_1 \cap Y_2} \} \right)$
  of $D_1$ with $(Y_1 \cap Y_2,|Y_1 \cap Y_2|) \subfan (Y_2,|Y_2|)$
  and a tropical structure that is compatible with the tropical structure on $D_2$.
  It is easily checked that $D_3$ is a refinement of $D_2$, too.
\end{proof}

\begin{definition}[Abstract tropical cycles] \label{defn-tropcycles}
  Let $((X,|X|),\omega_X)$ be an $n$-dimensional tropical polyhedral complex. Its equivalence class
  $\left[ \left( (X,|X|), \omega_X \right) \right]$ is called an \df{(abstract) tropical $n$-cycle}.
  The set of $n$-cycles is denoted by $Z_n$.
  Since the topological space $|X^{*}|$ of a tropical polyhedral complex $\left( (X,|X|), \omega_X \right)$
  is by definition invariant under refinements we define $\big|\left[ \left( (X,|X|), \omega_X \right) \right]\big|:=|X^{*}|$. Like in the affine case, an $n$-cycle $((X,|X|), \omega_X)$ is
  called an \df{(abstract) tropical variety} if $\omega_X(\sigma) \geq 0$ for all $\sigma \in X^{(n)}$.
  \medskip \\
  Let $C \in Z_n$ and $D \in Z_k$ be two tropical cycles. $D$ is called an
  \df{(abstract) tropical cycle in $C$} or a \df{subcycle of $C$} if there exists a representative
  $\left( ((Z,|Z|), \omega_{Z}), \{ \Psi_{\tau} \} \right)$ of $D$
  and a reduced representative $\left( ((X,|X|),\omega_{X}), \{ \Phi_{\sigma} \} \right)$ of $C$ such that
  \begin{enumerate}
    \item $(Z, |Z|) \subfan (X,|X|)$,
    \item the tropical structures on $Z$ and $X$ agree, i.e. for every $\tau \in Z$ the maps \linebreak
    ${\Psi_{\tau} \circ \Phi_{C_{Z,X}(\tau)}^{-1}}$ and $\Phi_{C_{Z,X}(\tau)} \circ \Psi_{\tau}^{-1}$
     are integer affine linear where defined.
  \end{enumerate}
  The set of tropical $k$-cycles in $C$ is denoted by $Z_k(C)$.
\end{definition}

\begin{remdef} \label{remdef-refiningembedcomplex}
  (a) Let $X$ be a finite set of rational polyhedra in $\RR^n$, $f \in
  \text{Hom}(\ZZ^n,\ZZ)$ a linear form and $b \in \RR$. Then let
  $$H_{f,b} := \big\{ \{x \in \RR^n| f(x) \leq b\}, \{x \in \RR^n| f(x)=b\}, \{x \in \RR^n| f(x) \geq b\}
  \big\}.$$ Like in the case of fans (cf. construction
  \ref{constr-refinement}) we can form sets $P :=$ \linebreak
  $\{\sigma \cap \sigma'| \sigma \in X, \sigma' \in H_{f,b} \}$ and
  $X \cap H_{f,b} := \{ \sigma \in P| \nexists \enspace \tau \in P^{(\dim(\sigma))} \text{ with } \tau \subsetneq \sigma  \}$.
  \medskip \\
  (b) Again let $X$ be a finite set of rational polyhedra in $\RR^n$. Let
  $\{f_i \leq b_i| i=1,\dots,N\}$ be all (integral) inequalities occurring in
  the description of all polyhedra in $X$. Then we can construct
  the set $X \cap H_{f_1,b_1} \cap \dots \cap H_{f_N,b_N}$. Note
  that for every collection of polyhedra $X$ this set $X \cap H_{f_1,b_1} \cap \dots \cap H_{f_N,b_N}$
  is a (usual) rational polyhedral complex (i.e. for every polyhedron $\tau \in X$ every
  face (in the usual sense) of $\sigma$ is contained in $X$ and
  the intersection of every two polyhedra in $X$ is a common face
  of each). Moreover note that the
  result is independent of the order of the $f_i$ and if $\{g_i \leq
  c_i| i=1,\dots,M\}$ is a different set of inequalities describing the
  polyhedra in $X$ then $X \cap H_{f_1,b_1} \cap \dots \cap H_{f_N,b_N}$
  and $X \cap H_{g_1,c_1} \cap \dots \cap H_{g_M,c_M}$ have a common
  refinement, namely $X \cap H_{f_1,b_1} \cap \dots \cap H_{f_N,b_N} \cap H_{g_1,c_1} \cap \dots \cap
  H_{g_M,c_M}$.
\end{remdef}

\begin{construction}[Sums of tropical cycles]\label{constr-sumsofcycles}
  Let $C \in Z_n$ be a tropical cycle. Like in the affine case the set
  of tropical $k$-cycles in $C$ can be made into an abelian group by defining the sum of two such $k$-cycles
  as follows: Let $D_1$ and $D_2 \in Z_k(C)$ be the two cycles whose sum we want to construct.
  By definition there are reduced representatives
  $\left( ((X_1,|X_1|),\omega_{X_1}), \{\Phi_{\tau}^{X_1}\} \right)$ and
  $\left( ((X_2,|X_2|),\omega_{X_2}), \{\Phi_{\tau}^{X_2}\} \right)$
  of $C$ and reduced representatives $\left( ((Y,|Y|),\omega_{Y}), \{\Phi_{\tau}^{Y}\} \right)$ of
  $D_1$ and $\left( ((Z,|Z|), \omega_{Z}), \{\Phi_{\tau}^{Z}\} \right)$ of $D_2$ such that
  $(Y,|Y|) \subfan (X_1,|X_1|)$ and the tropical structures on $Y$ and $X_1$ agree and
  $(Z,|Z|) \subfan (X_2,|X_2|)$ and the tropical
  structures on $Z$ and $X_2$ agree. As ``$\sim$'' is an equivalence relation
  there is a common refinement $\left( ((X,|X|, \{ \varphi_{\tau} \}),\omega_{X}), \{\Phi_{\tau}^{X}\} \right)$ of $X_1$
  and $X_2$ which we may assume to be reduced.
  Applying construction \ref{constr-constrrefinemofwpolycompl} to $Y$ and $X$
  we obtain the tropical polyhedral complex $\left( (({Y \cap X}, {|Y \cap X|}), \omega_{Y \cap X}), \{\Phi_{\tau}^{Y \cap X}\} \right)$
  which is a refinement of $Y$, has a tropical structure that is compatible with the tropical structure on $X$
  and fulfils ${(Y \cap X, |Y \cap X|)} \subfan (X,|X|)$. If we further apply construction
  \ref{constr-constrrefinemofwpolycompl} to $Z$ and $X$ we get a
  refinement $\left( ((Z \cap X, |Z \cap X|), \omega_{Z \cap X}), \{\Phi_{\tau}^{Z \cap X}\} \right)$ of $Z$
  with analogous properties. Now fix some polyhedron $\sigma \in X$ and let $\tau_1, \dots, \tau_r \in Y \cap X$
  and $\tau_{r+1}, \dots, \tau_s \in Z \cap X$ be all polyhedra of $Y \cap X$ and $Z \cap X$ respectively
  that are contained in $\sigma$. Note that property (a) of definition
  \ref{defn-tropcycles} implies that for all $i=1, \dots, r$ the image $\varphi_{\sigma}(\tau_i)$ is
  a rational polyhedron in $\RR^{n_{\sigma}}$. Like in remark and definition \ref{remdef-refiningembedcomplex}
  let $\{f_i \leq b_i| i=1,\dots,N\}$ be the set of all integral inequalities occurring in the
  description of all polyhedra $\varphi_{\sigma}(\tau_i), i=1, \dots, s$
  and let $R_{Y \cap X}^{\sigma} := \{\varphi_{\sigma}(\tau_i)| i=1, \dots, r \} \cap H_{f_1,b_1} \cap
  \dots \cap H_{f_N,b_N}$ and $R_{Z \cap X}^{\sigma} := \linebreak \{\varphi_{\sigma}(\tau_i)| i=r+1, \dots, s \} \cap H_{f_1,b_1} \cap \dots
  \cap H_{f_N,b_N}$. Then $P_{Y \cap X}^{\sigma} := \{ \varphi_{\sigma}^{-1}(\tau)| \tau \in
  R_{Y \cap X}^{\sigma} \}$ and $P_{Z \cap X}^{\sigma} := \{ \varphi_{\sigma}^{-1}(\tau)| \tau \in
  R_{Z \cap X}^{\sigma} \}$ are a kind of local refinement of $Y \cap X$ and $Z \cap
  X$ respectively, but taking the union over all maximal polyhedra
  $\sigma \in X^{(n)}$ does in general not lead to global refinements as there may
  be overlaps between polyhedra coming from different $\sigma$.
  We resolve this as follows:
  For $\sigma \in X^{(n)}$, $\tau \in \bigcup_{i=0}^{n-1} X^{(i)}$ let \linebreak $P_{Y,\tau}^{\sigma}
  := \{ \varrho \in P_{Y \cap X}^{\sigma}| \tau \text{ is the inclusion-minimal polyhedron of } X \text{ containing } \varrho\}$
  and \linebreak $P_{Y,n} := \bigcup_{\sigma \in X^{(n)}} \{ \varrho \in P_{Y \cap X}^{\sigma}
  | \nexists \widetilde{\tau} \in X^{(n-1)}: \varrho \subseteq \widetilde{\tau}
  \}$. Analogously for $P_{Z,\tau}^{\sigma}$ and $P_{Z,n}$. \linebreak Then let
  $\widetilde{Y}:= P_{Y,n} \cup \left( \bigcup_{\tau \in X^{(i)}: i<n} \{ \bigcap_{\sigma \in X^{(n)}: \tau \subseteq \sigma}
  \tau_\sigma| \tau_\sigma \in P_{Y,\tau}^{\sigma} \} \right)$
  and $\widetilde{Z}:= \linebreak P_{Z,n} \cup \left( \bigcup_{\tau \in X^{(i)}: i<n} \{\bigcap_{\sigma \in X^{(n)}: \tau \subseteq \sigma}
  \tau_\sigma| \tau_\sigma \in P_{Z,\tau}^{\sigma}\} \right)$. Moreover for every $\tau \in \widetilde{Y} \cup \widetilde{Z}$ choose some
  $\sigma \in X^{(n)}$ with $\tau \subseteq \sigma$ and let $\psi_{\tau}:=\varphi_{\sigma}|_{\tau}$.
  Note that by construction $(\widetilde{Y},|Y \cap X|)$ and $(\widetilde{Z},|Z \cap X|)$ with
  structure maps $\psi_{\tau}, \tau \in \widetilde{X}$ or $\tau \in \widetilde{Z}$ respectively and weight functions
  $\omega_{\widetilde{Y}}$ and $\omega_{\widetilde{Z}}$ induced by $Y \cap X$ and $Z \cap X$ are refinements of
  $Y \cap X$ and $Z \cap X$ (we need here that $R_{Y \cap X}^{\sigma}$ and $R_{Z \cap X}^{\sigma}$ were usual polyhedral
  complexes in $\RR^{n_{\sigma}}$). Thus we can endow them with the tropical structures inherited from $Y \cap X$ and $Z \cap X$
  respectively (cf. construction \ref{constr-refinementtropcomplex}).
  As $(\widetilde{X} \cup \widetilde{Y},|Y \cap X| \cup |Z \cap X|)$ is a polyhedral complex
  \begin{center}
  \input{pic/Subdiv.pstex_t}\\
  An illustration of the process described in construction
  \ref{constr-sumsofcycles}.
  \end{center}
  now, we can form $$((P,|P|),\omega_P) := ((\widetilde{X} \cup \widetilde{Y},|Y \cap X| \cup |Z \cap
  X|), \omega_P),$$ where $\omega_{P}(\sigma) := \omega_{\widetilde{Y}}(\sigma) + \omega_{\widetilde{Z}}(\sigma)$
  for all $\sigma \in P^{(k)}$ (we set $\omega_{\Box}(\sigma) := 0$ for ${\sigma \notin \Box},
  \linebreak \Box \in \{\widetilde{Y}, \widetilde{Z} \}$).
  Recall that the tropical structures on $\widetilde{Y}$ and $\widetilde{Z}$ are inherited
  from $Y \cap X$ and $Z \cap X$ and are thus compatible with the
  tropical structure on $X$. Thus $\Phi_{\sigma}^{X}(S_{\sigma}^P) \subseteq |F_{\sigma}^{X}|$
  with weights induced from $P$ is an open fan (the corresponding complete tropical fan is just the sum of the fans
  coming from $\widetilde{Y}$ and $\widetilde{Z}$). Thus we can set
  $\widetilde{\Phi}_{\sigma} := \Phi_{\sigma}^{X} |_{S_{\sigma}^P}:
  S_{\sigma}^P \stackrel{\sim}{\rightarrow}
  \Phi_{\sigma}^{X}(S_{\sigma}^P)$ and can hence define the sum $D_1+D_2$ to be
  $$ D_1 + D_2 := \left[ \left( ((P,|P|),\omega_P), \{ \widetilde{\Phi}_{\sigma} \} \right) \right].$$
  Note that the class $[( ((P,|P|),\omega_P), \{ \widetilde{\Phi}_{\sigma} \})]$
  is independent of the choices we made, i.e. the sum $D_1 + D_2$ is well-defined.
\end{construction}

\begin{lemma} \label{lemma-cyclegroup}
  Let $C \in Z_n$ be a tropical cycle. The set $Z_k(C)$ together with the operation
  ``+'' from construction \ref{constr-sumsofcycles} forms an abelian group.
\end{lemma}
\begin{proof}
  The class of the empty complex $0=[\emptyset]$ is the neutral
  element of this operation and $[((Y,|Y|),-\omega_Y)]$ is the inverse element of
  $[((Y,|Y|),\omega_Y)] \in Z_k(C)$.
\end{proof}

  \section{Cartier divisors and their associated Weil divisors} \label{sec-cartierdivisors}

\begin{definition}[Rational functions and Cartier divisors] \label{def-rationalfunctionsandcartierdivisors}
    Let $C$ be an abstract $k$-cycle and let $U$ be an open set in $|C|$. A \df{(non-zero) rational function on $U$} is a continuous function $\varphi : U \rightarrow \RR$ such that there exists a representative $(((X, |X|, \{m_\sigma\}_{\sigma \in X}), \omega_X), \{M_\sigma\}_{\sigma \in X})$ of $C$ such that for each face $\sigma \in X$ the map $\varphi \circ m_\sigma^{-1}$ is locally integer affine linear (where defined). The \df{set of all non-zero rational functions on $U$} is denoted by $\calK^*_C(U)$ or just $\calK^*(U)$. \\
    If additionally for each face $\sigma \in X$ the map $\varphi \circ M_\sigma^{-1}$ is locally integer affine linear (where defined), $\varphi$ is called \df{regular invertible}. The \df{set of all regular invertible functions on $U$} is denoted by $\calO^*_C(U)$ or just $\calO^*(U)$. \\
    A \df{representative of a Cartier divisor on $C$} is a finite set $\{(U_1, \varphi_1), \ldots , (U_l, \varphi_l)\}$, where $\{U_i\}$ is an open covering of $|C|$ and $\varphi_i \in \calK^*(U_i)$ are rational functions on $U_i$ that only differ in regular invertible functions on the overlaps, in other words, for all $i \neq j$ we have $\varphi_i|_{U_i \cap U_j} - \varphi_j|_{U_i \cap U_j} \in \calO^*(U_i \cap U_j)$. \\
    We define the \df{sum} of two representatives by $\{(U_i, \varphi_i)\} + \{(V_j, \psi_j)\} = \{(U_i \cap V_j, \varphi_i + \psi_j)\}$, which obviously fulfills again the condition on the overlaps. \\
    We call two representatives $\{(U_i, \varphi_i)\}, \{(V_j, \psi_j)\}$ \df{equivalent} if $\varphi_i - \psi_j$ is regular invertible (where defined) for all $i,j$, i.e. $\{(U_i, \varphi_i)\} - \{(V_j, \psi_j)\} = \{(W_k, \gamma_k)\}$ with $\gamma_k \in \calO^*(W_k)$. Obviously, ``+'' induces a group structure on the set of equivalence classes of representatives with the neutral element $\{(|C|, c_0)\}$, where $c_0$ is the constant zero function. This group is denoted by $\Div (C)$ and its elements are called \df{Cartier divisors on $C$}.
\end{definition}

\begin{example} \label{example-cartierdivisor}
    \begin{figure}[t]
      \input{pic/ExCartDiv.pstex_t}\\
      The Cartier divisor $\varphi$ defined in example \ref{example-cartierdivisor}.
    \end{figure}
    Let us give an example of a Cartier divisor which is not globally defined by a rational function: As abstract cycle $C$ we take the elliptic curve $[X_2]$ from example \ref{example-abstractcycles} (the brackets resemble the fact that, to be precise, we take the equivalence class of the polyhedral complex $X_2$ with respect to refinements). By $\alpha_1, \alpha_2$ we denote the two vertices in $X_2$. W.l.o.g. we can assume that the maps $M_{\alpha_i}$ map the points $\alpha_i$ exactly to $0 \in \RR$. Of course, the stars $S_{\alpha_1}, S_{\alpha_2}$ cover our whole space $|C| = |X_2|$. So we can define the Cartier divisor $\varphi := [\{(S_{\alpha_1}, \psi_1), (S_{\alpha_2}, \psi_2)\}]$, where $\psi_1 := \max(0,x) \circ M_{\alpha_1}$ and $\psi_2 := c_0 \circ M_{\alpha_2}$ with $c_0$ the constant zero function. Let us check the condition on the overlaps: On one open half of our curve the two functions coincide, whereas on the other open half they differ by a linear function. So we constructed an Cartier divisor which can not be globally defined by one rational function (as $\psi_1$ can not be completed to a continuous function on $|C|$).
\end{example}

\begin{remark}[Restrictions to subcycles] \label{remark-restrictiontosubcycles}
Note that, as in the affine case (see remark \ref{remark-zerofunction}), we can restrict a non-zero rational function $\varphi \in \calK^*_C(U)$ to an arbitrary subcycle $D \subseteq C$, i.e. $\varphi|_{U \cap |D|} \in \calK^*_D(U \cap |D|)$. It is also true that a regular invertible function $\varphi \in \calO^*_C(U)$ restricted to $D$ is again regular invertible, i.e. $\varphi|_{U \cap |D|} \in \calO^*_D(U \cap |D|)$. Hence we can also restrict a Cartier divisor $[\{(U_i, \varphi_i)\}] \in \Div(C)$ to $D$ by setting $[\{(U_i, \varphi_i)\}]\,|_D := [\{(U_i \cap |D|, \varphi_i|_{U_i \cap |D|})\}] \in \Div(D)$. Let us also stress again that we still require our objects to be defined everywhere (on a given open subset $U$). This causes problems like for example in remark \ref{remark-pushforwardratequiv}.
\end{remark}

\begin{construction}[Intersection products] \label{constr-intersectionproducts}
    Let $C$ be an abstract $k$-cycle and $\varphi = [\{(U_i, \varphi_i)\}] \linebreak \in \Div(C)$ a Cartier divisor on $C$. By  definition \ref{def-rationalfunctionsandcartierdivisors} and lemma \ref{lemma-equivrelontroppolycomplex}, there exists a representative $(((X, |X|, \{m_\sigma\}_{\sigma \in X}), \omega_X), \{M_\sigma\}_{\sigma \in X})$ of $C$ such that for all $i$ and $\sigma \in X$ the map $\varphi_i \circ m_\sigma^{-1}$ is locally integer affine linear (where defined). We can also assume that $X = X^*$, as our functions are defined on $|C| = |X^*|$ at the most. We would like to define the intersection product $\varphi \cdot C$ to be
  $$
    \bigg[\bigg(\Big(\big(Y, |Y|, \{m_\sigma\}_{\sigma \in Y}\big), \omega_{X,\varphi}\Big), \big\{M_\sigma|_{S^Y_\sigma} : S^Y_\sigma \rightarrow |F^Y_\sigma|\big\}_{\sigma \in Y}\bigg)\bigg],
  $$
  where
  $$
    Y := \bigcup_{i=0}^{k-1} X^{(i)}, \;
    |Y| := \bigcup_{\sigma \in Y} \sigma, \;
    S^Y_\sigma = \bigcup_{\substack{\sigma' \in Y \\ \sigma \subseteq \sigma'}} (\sigma')^{ri}, \;
    F^Y_\sigma := \bigcup_{i=0}^{k-1} F_\sigma^{(i)}
  $$
    and $\omega_{X,\varphi}$ is an appropriate weight function. So it remains to construct $\omega_{X,\varphi}(\tau)$ for $\tau \in X^{(k-1)}$. \\
    First, we do this pointwise, i.e. we construct $\omega_{X,\varphi}(p)$ for $p \in (\tau)^{ri}$. Given a $p \in (\tau)^{ri}$, we pick an $i$ with $p \in U_i$. Let $V$ be the connected component of $M_\tau(U_i \cap S_\tau)$ containing $M_\tau(p)$. Then the function $\varphi_i \circ M_\tau^{-1}|_V$ can be uniquely extended to a rational function $\tilde{\varphi}_i \in \calK^*([(\tilde{F}_\tau, \omega_{\tilde{F}_\tau})])$, where $(\tilde{F}_\tau, \omega_{\tilde{F}_\tau})$ is the tropical fan generated by the open fan $(F_\tau, \omega_{F_\tau})$. So, in the affine case, we can compute $\omega_{\tilde{F}_\tau, \tilde{\varphi}_i}(\RR \cdot M_\tau(\tau))$ (see construction \ref{constr-associatedweildivisor} and definition \ref{def-associatedweildivisor}) and define $\omega_{X,\varphi}(p) := \omega_{\tilde{F}_\tau, \tilde{\varphi}_i}(\RR \cdot M_\tau(\tau))$. \\
    This definition is well-defined, namely if we pick another $j$ with $p \in U_j$ and denote by $V'$ the connected component of $M_\tau(U_j \cap S_\tau)$ containing $M_\tau(p)$, we know by definition of a Cartier divisor that $\varphi_i \circ M_\tau^{-1}|_{V \cap V'} - \varphi_j \circ M_\tau^{-1}|_{V \cap V'}$ is affine linear, hence $\tilde{\varphi}_i - \tilde{\varphi}_j$ is affine linear. By remark \ref{rem-affinelinearfunctionsandsums} we get $\omega_{\tilde{F}_\tau, \tilde{\varphi}_i}(\RR \cdot M_\tau(\tau)) = \omega_{\tilde{F}_\tau, \tilde{\varphi}_j}(\RR \cdot M_\tau(\tau))$. \\
    The same argument shows that our definition does not depend on the choice of a representative $\{(U_i, \varphi_i)\}$ of $\varphi$. \\
    But as $(\tau)^{ri}$ is connected, the continuous function $\omega_{X,\varphi} : (\tau)^{ri} \rightarrow \ZZ$ must be constant. Hence, we define $\omega_{X,\varphi}(\tau) := \omega_{X,\varphi}(p)$ for some $p \in (\tau)^{ri}$. With this weight function
    $$
    \bigg(\Big(\big(Y, |Y|, \{m_\sigma\}_{\sigma \in Y}\big), \omega_{X,\varphi}\Big), \{M_\sigma|_{S^Y_\sigma}\}_{\sigma \in Y}\bigg)
  $$
    is a tropical polyhedral complex. \\
    Let us now check if the equivalence class of this complex is independent of the choice of representatives of $C$. Let therefore $(((X', |X'|, \{m_{\sigma'}\}_{\sigma' \in X'}), \omega_{X'}), \{M_{\sigma'}\}_{\sigma' \in X'})$ be a refinement of $(((X, |X|, \{m_\sigma\}_{\sigma \in X}), \omega_X), \{M_\sigma\}_{\sigma \in X})$ (we can again assume $X'=X'^*$). Then, for each $\sigma' \in X'$, the map $M_{C_{X',X}(\sigma')} \circ M_{\sigma'}^{-1}$ embeds $F_{\sigma'}$ into a refinement of $F_{C_{X',X}(\sigma')}$. Applying the affine statement here (see remark \ref{rem-welldefined}), we deduce that for each $\tau' \in X'^{(k-1)}$ it holds $\omega_{X', \varphi}(\tau') = 0$ (if $\dim C_{X',X}(\tau') = k$) or $\omega_{X', \varphi}(\tau')) = \omega_{X, \varphi}(C_{X',X}(\tau'))$ (if $\dim C_{X',X}(\tau') = k-1$).
\end{construction}

\begin{definition}[Intersection products] \label{def-intersectionproducts}
  Let $C$ be an abstract $k$-cycle and $\varphi = [\{(U_i, \varphi_i)\}] \in \Div(C)$ a Cartier divisor on $C$. Let furthermore $(((X, |X|, \{m_\sigma\}_{\sigma \in X}), \omega_X), \{M_\sigma\}_{\sigma \in X})$ be a representative of $C$ such that $|X| = |C|$ and for all $i$ and $\sigma \in X$ the map $\varphi_i \circ m_\sigma^{-1}$ is locally integer affine linear (where defined). The \df{associated Weil divisor} $\divisor(\varphi) = \varphi \cdot C$ is defined to be
  $$
    \bigg[\bigg(\Big(\big(Y := \bigcup_{i=0}^{k-1} X^{(i)}, \bigcup_{\sigma \in Y} \sigma, \{m_\sigma\}_{\sigma \in Y}\big), \omega_{X,\varphi}\Big), \{M_\sigma|_{S^Y_\sigma}\}_{\sigma \in Y}\bigg)\bigg] \in Z_{k-1}(C),
  $$
  where $S^Y_\sigma = \bigcup_{\substack{\sigma' \in Y \\ \sigma \subseteq \sigma'}} (\sigma')^{ri}$ and $\omega_{X,\varphi}$ is the weight function constructed in construction \ref{constr-intersectionproducts}. \\
  Let $D \in Z_l(C)$ be an arbitrary subcycle of $C$ of dimension $l$. We define the \df{intersection product of $\varphi$ with $D$} to be $\varphi \cdot D:= \varphi|_D \cdot D \in Z_{l-1}(C)$.
\end{definition}

\begin{example}
    Let us compute the Weil divisor associated to our Cartier divisor $\varphi$ on the elliptic curve $C$ constructed in example \ref{example-cartierdivisor}. In fact, there is nothing to compute: One can see immediately from the picture that $\divisor(\varphi)$ is just the vertex $\alpha_1$ with multiplicity 1 (the multiplicity of $\alpha_2$ is 0 as in order to compute it, one has to use the constant function $\psi_2$). Let us stress that this single point can not be obtained as the Weil divisor of a (global) rational function, as all such divisors must have ``degree 0'' (this is defined precisely and proven in remark \ref{redundantremark} and lemma \ref{lemma-degreezero}).
\end{example}

\begin{proposition}[Commutativity]
\label{prop-commutativity}
    Let $\varphi, \psi \in \Div(C)$ be two Cartier divisors on $C$. Then $\psi \cdot (\varphi \cdot C) = \varphi \cdot (\psi \cdot C)$.
\end{proposition}

\begin{proof}
    Say $\varphi = [\{(U_i, \varphi_i)\}]$ and $\psi = [\{(V_j, \psi_j)\}]$. Using lemma \ref{lemma-equivrelontroppolycomplex} we find a representative $(((X, |X|, \{m_\sigma\}_{\sigma \in X}), \omega_X), \{M_\sigma\}_{\sigma \in X})$ of $C$ such that $|X| = |C|$ and for all $i,j$ and $\sigma \in X$ the maps $\varphi_i \circ m_\sigma^{-1}$ and $\psi_j \circ m_\sigma^{-1}$ are locally integer affine linear (where defined). For $\theta \in X^{(k-2)}$, $p \in (\theta)^{ri}$ and $i,j$ with $p \in U_i \cap V_j$ we get (using notations from construction \ref{constr-intersectionproducts}) $\omega_{X, \varphi, \psi}(\theta) = \omega_{X, \varphi, \psi}(p) = \omega_{\tilde{F}_\theta, \tilde{\varphi}_i, \tilde{\psi}_j}(\RR \cdot M_\theta(\theta))$ and similarily $\omega_{X, \psi, \varphi}(\theta) = \omega_{\tilde{F}_\theta, \tilde{\psi}_j, \tilde{\varphi}_i}(\RR \cdot M_\theta(\theta))$. Using the corresponding statement in the affine case now (see proposition \ref{prop-balancingconditionandcommutativity} (b)), we deduce that the two weight functions are equal which proves the claim.
\end{proof}

  \section{Push-forward of tropical cycles and pull-back of Cartier divisors} \label{sec-gluedpushpull}

\begin{definition}[Morphisms of tropical cycles] \label{def-morphoftropcycle}
  Let $C \in Z_n$ and $D \in Z_m$ be two tropical cycles. A \df{morphism} $f: C \rightarrow
  D$ of tropical cycles is a continuous map $f: |C| \rightarrow |D|$
  with the following property:
  There exist reduced representatives $\left(((X,|X|), \omega_X), \{\Phi_{\sigma}\} \right)$
  of $C$ and $\left(((Y,|Y|),\omega_Y),\{\Psi_{\tau}\} \right)$ of $D$ such that
  \begin{enumerate}
    \item \label{pma} for every polyhedron $\sigma \in X$ there exists a polyhedron $\widetilde{\sigma} \in Y$
    with $f(\sigma) \subseteq \widetilde{\sigma}$,
    \item for every pair $\sigma, \widetilde{\sigma}$ from \ref{pma} the map
    $\Psi_{\widetilde{\sigma}} \circ f \circ \Phi_{\sigma}^{-1}: |F_{\sigma}^{X}| \rightarrow |F_{\widetilde{\sigma}}^{Y}|$
    induces a morphism of fans $\widetilde{F}_{\sigma}^{X} \rightarrow
    \widetilde{F}_{\widetilde{\sigma}}^{Y}$ (cf. definition \ref{morphism}), where $\widetilde{F}_{\sigma}^{X}$ and
    $\widetilde{F}_{\widetilde{\sigma}}^{Y}$ are the tropical fans associated to $F_{\sigma}^{X}$ and $F_{\widetilde{\sigma}}^{Y}$
    respectively (cf. definition \ref{defn-openfan}).
  \end{enumerate}
\end{definition}

First of all we want to show that the restriction of a morphism to
a subcycle is again a morphism:

\begin{lemma} \label{lemma-restrtosubcycles}
  Let $C \in Z_n$ and $D \in Z_m$ be two cycles, $f: C \rightarrow D$
  a morphism and $E \in Z_k(C)$ a subcycle of $C$. Then the map
  $f|_{|E|} : |E| \rightarrow |D|$ induces a morphism of tropical
  cycles $f|_{E}: E \rightarrow D$.
\end{lemma}
\begin{proof}
  By definition of morphism there exist reduced representatives $((X_1,|X_1|),
  \omega_{X_1})$ of $C$ and $((Y,|Y|),\omega_{Y})$ of $D$ such that
  properties (a) and (b) in definition \ref{def-morphoftropcycle}
  are fulfilled. By definition of subcycle there exist reduced
  representatives $((Z_1,|Z_1|),\omega_{Z_1})$ of $E$ and
  $((X_2,|X_2|),\omega_{X_2})$ of $C$ such that properties (a) and
  (b) in definition \ref{defn-tropcycles} are fulfilled, i.e.
  such that $(Z_1,|Z_1|) \subfan (X_2,|X_2|)$ and the
  tropical structures on $Z_1$ and $X_2$ agree.
  As ``$\sim$'' is an equivalence relation there exists a common refinement
  $((X,|X|),\omega_{X})$ of $((X_1,|X_1|),\omega_{X_1})$ and
  $((X_2,|X_2|),\omega_{X_2})$ which we may assume to be reduced.
  Applying construction \ref{constr-constrrefinemofwpolycompl}
  to $Z_1$ and $X$ we obtain a refinement $((Z,|Z|),\omega_Z) := \linebreak
  {((Z_1 \cap X,|Z_1 \cap X|),\omega_{Z_1 \cap
  X})}$ of $((Z_1,|Z_1|), \omega_{Z_1})$ such that $(Z,|Z|) \subfan
  (X,|X|)$ and the tropical structures on $Z$ and $X$ agree.
  Thus properties (a) and (b) of definition \ref{def-morphoftropcycle}
  are fulfilled by $Z$ and $Y$ and the restricted map
  $f|_{|E|} : |E| \rightarrow |D|$ gives us a morphism $f|_{E}: E \rightarrow D$.
\end{proof}

If we are given a morphism and a tropical cycle the following
construction shows how to build the push-forward cycle of the
given one along our morphism:

\begin{construction}[Push-forward of tropical cycles] \label{constr-pushfoftropcycles}
  Let $C \in Z_n$ and $D \in Z_m$ be two cycles and let $f: C \rightarrow D$
  be a morphism. Let $\left(((X,|X|, \{\varphi_{\sigma}\}), \omega_X), \{\Phi_{\sigma}\} \right)$
  and $\left(((Y,|Y|,\{\psi_{\sigma}\}),\omega_Y),\{\Psi_{\tau}\} \right)$ be
  representatives of $C$ and $D$ fulfilling properties (a) and (b)
  of definition \ref{def-morphoftropcycle}. Consider the collection of polyhedra $$Z := \{
  f(\sigma) | \sigma \in X \text{ contained in a maximal
  polyhedron of } X \text{ on which } f \text{ is injective} \}.$$
  In general $Z$ is not a polyhedral complex. We resolve this by
  subdividing the polyhedra in $Z$ and refining $X$ accordingly:\\
  Fix some polyhedron $\widetilde{\sigma} \in
  Y^{(m)}$ and let $\tau_1,\dots,\tau_r \in Z$ be all polyhedra that are contained
  in $\widetilde{\sigma}$. Property (b) of definition
  \ref{def-morphoftropcycle} implies that $\{
  \psi_{\widetilde{\sigma}}(\tau_i) | i=1,\dots,r\}$ is a set of
  rational polyhedra in $\RR^{n_{\widetilde{\sigma}}}$. Like in
  remark and definition \ref{remdef-refiningembedcomplex} let
  $\{ g_i(x) \leq b_i| i=1,\dots,N\}$, $g_i \in
  \text{Hom}(\ZZ^{n_{\widetilde{\sigma}}},\ZZ)$, $b_i \in \RR$
  be all inequalities occurring in the description
  of all polyhedra in $\{\psi_{\widetilde{\sigma}}(\tau_i) |
  i=1,\dots,r\}$ and let
  \begin{eqnarray*}
  R_{\widetilde{\sigma}} & := & \{\psi_{\widetilde{\sigma}}(\tau_i) | i=1,\dots,r\} \cap H_{G_{1},b_1} \cap \dots \cap
  H_{G_{N},b_N},\\
  P_{\widetilde{\sigma}} & := & \{ \psi_{\widetilde{\sigma}}^{-1}(\tau)| \tau \in
  R_{\sigma_i}\}.
  \end{eqnarray*}
  Like in construction \ref{constr-sumsofcycles}
  $P_{\widetilde{\sigma}}$ can be seen as a kind of local refinement of
  $Z$. But here again taking the union over all maximal polyhedra
  $\widetilde{\sigma} \in Y^{(m)}$ does in general not lead to a global refinement as there may
  be overlaps between polyhedra coming from different $\widetilde{\sigma}$. We
  fix this as follows (cf. \ref{constr-sumsofcycles}):
  For $\widetilde{\sigma} \in Y^{(m)}$ and $\widetilde{\tau} \in \bigcup_{i=0}^{m-1} Y^{(i)}$ let \linebreak $P_{Z,\widetilde{\tau}}^{\widetilde{\sigma}} :=
  \{ \varrho \in P_{\widetilde{\sigma}} | \widetilde{\tau} \text{ is the inclusion minimal polyhedron of }
  Y \text{ containing } \varrho \} $
  and \linebreak $P_{Z,m} := \bigcup_{\widetilde{\sigma} \in Y^{(m)}}
  \{ \varrho \in P_{\widetilde{\sigma}} | \nexists \widetilde{\tau} \in Y^{(m-1)}:
  \varrho \subseteq \widetilde{\tau} \}$.
  Then $\widetilde{Z}:= \linebreak P_{Z,m} \cup \left( \bigcup_{\widetilde{\tau} \in Y^{(i)}: i<m}
  \{\bigcap_{\widetilde{\sigma} \in Y^{(m)}: \widetilde{\tau} \subseteq \widetilde{\sigma}}
  \tau_{\widetilde{\sigma}}| \tau_{\widetilde{\sigma}} \in P_{Z,\widetilde{\tau}}^{\widetilde{\sigma}}\} \right)$ is
  the set of polyhedra (without any overlaps now) that shall
  induce our wanted refinement of $X$: Let $T := \linebreak \{ \sigma \in X^{(n)}| f \text{ is injective on } \sigma \}$,
  $Q_0 := \{\tau \in X| \nexists \sigma \in T: \tau \subseteq \sigma\}$ and
  $Q_1 := \linebreak \left( \bigcup_{\sigma \in T} \{ (f|_{\sigma})^{-1} (\tau) |
  \tau \in \widetilde{Z}, \tau \subseteq f(\sigma) \} \right)$.
  Then define $\widetilde{X}:= Q_0 \cup Q_1$.\\
  Let $\tau \in Q_1$ and choose $\sigma \in T$ with
  $\tau \subseteq \sigma$. Property (b) of definition \ref{def-morphoftropcycle} implies that
  $\psi_{\widetilde{\sigma}} \circ f \circ \varphi_{\sigma}^{-1}$ is integer affine linear where
  defined. Hence $\varphi_{\sigma}(\tau)$ is a rational polyhedron
  in $\RR^{n_{\sigma}}$. Denote by $H_{\sigma,\tau}$ the smallest
  affine subspace of $\RR^{n_{\sigma}}$ containing
  $\varphi_{\sigma}(\tau)$. We can consider $\varrho_{\tau} :=
  \varphi_{\sigma}|_{\tau}$ to be a map $\varrho_{\tau} : \tau \rightarrow
  H_{\sigma,\tau} \cong \RR^{n_\tau}$. Note that by construction
  $(\widetilde{X},|X|,\{\varrho_{\tau}\})$ is a polyhedral
  complex. We endow it with the weight function $\omega_{\widetilde{X}}$
  and tropical structure $\{ \Phi_{\tau}^{\widetilde{X}} \}$ induced by $X$. Now we are able to define $$f_{*}X := \{
  f(\sigma) | \sigma \in \widetilde{X} \text{ contained in a maximal
  polyhedron of } \widetilde{X} \text{ on which } f \text{ is injective}
  \}$$ and $|f_{*}X| := \bigcup_{\tau \in f_{*}X} \tau$. For every polyhedron $\tau \in f_{*}X$
  let $\sigma_{\tau} \in Y$ be the inclusion-minimal polyhedron containing $\tau$.
  Then define $\vartheta_{\tau} := \psi_{\sigma_\tau} |_{\tau}: \tau \rightarrow H_{\sigma_\tau,\tau}
  \cong \RR^{n_{\tau}}$, where $H_{\sigma_\tau,\tau} \subseteq \RR^{n_{\sigma_\tau}}$ is the smallest affine subspace
  containing the rational polyhedron $\psi_{\sigma_\tau}(\tau) \in \widetilde{Z}$. Note that this makes
  $(f_{*}X,|f_{*}X|,\{\vartheta_{\tau}\})$ into a polyhedral
  complex. Moreover note that property (b) of definition
  \ref{def-morphoftropcycle} still holds for $\widetilde{X}$ and
  $Y$. Hence we can assign weights and tropical fans to
  $f_{*}X$ as follows: Let $\sigma \in f_{*}X$, let $\widetilde{\sigma} \in Y$ be the inclusion-minimal
  polyhedron containing it and let
  $\tau_1, \dots, \tau_r \in \widetilde{X}$ be all polyhedra with
  $f(\tau_i)=\sigma$ that are contained in a maximal polyhedron
  of $\widetilde{X}$ on which $f$ is injective. Then let
  $\Psi_{\widetilde{\sigma}}(S_{\widetilde{\sigma}})=F_{\widetilde{\sigma}}^{Y}$
  and $\Phi_{\tau_i}^{\widetilde{X}}(S_{\tau_i})=F_{\tau_i}^{\widetilde{X}}$ respectively
  be the corresponding open fans and $\widetilde{F}_{\widetilde{\sigma}}^{Y}$,
  $\widetilde{F}_{\tau_i}^{\widetilde{X}}$
  be the associated tropical fans. Property (b) of definition
  \ref{def-morphoftropcycle} implies that
  $f_{*}\widetilde{F}_{\tau_i}^{\widetilde{X}} \subseteq |\widetilde{F}_{\widetilde{\sigma}}^{Y}|$ is again a
  tropical fan (note that we do not need to refine $\widetilde{F}_{\tau_i}^{\widetilde{X}}$
  to construct this push-forward). Thus we can define
  $$\left(\widetilde{F}_{\sigma}^{f_{*}X}, \omega_{\widetilde{F}_{\sigma}^{f_{*}X}}\right):= \left(\bigcup_{i=1}^r
  f_{*}\widetilde{F}_{\tau_i}^{\widetilde{X}}, \sum_{i=1}^r \omega_{f_{*}\widetilde{F}_{\tau_i}^{\widetilde{X}}} \right)
  \quad \text{and} \quad F_{\sigma}^{f_{*}X} :=
  \widetilde{F}_{\sigma}^{f_{*}X} \cap \Psi_{\widetilde{\sigma}}(S_{\sigma})$$
  (here again we assume that $\omega_{f_{*}\widetilde{F}_{\tau_i}^{\widetilde{X}}}(\tau)=0$ if
  $\tau \notin f_{*}\widetilde{F}_{\tau_i}^{\widetilde{X}}$). Moreover we define
  $$ \Theta_{\sigma} := \Psi_{\widetilde{\sigma}}|_{S_\sigma}: S_{\sigma} \rightarrow |F_{\sigma}^{f_{*}X}|.$$
  Then the map $\Theta_{\sigma}$, $\sigma \in f_{*}X$ is 1:1 on polyhedra and we can endow the maximal polyhedra of $f_{*}X$ with
  weights $\omega_{f_{*}X}(\cdot)$ coming from $F_{\sigma}^{f_{*}X}$ in this way. These weights are obviously
  well-defined by property (c) of the tropical polyhedral complex $Y$ (cf. definition \ref{defn-tropcomplex})
  and the maps $\Theta_{\sigma}$ for different $\sigma \in f_{*}X$ are obviously
  compatible. Hence we can define
  $$f_{*}C := \bigg[ \big( \big( (f_{*}X,|f_{*}X|,\{ \vartheta_\tau \}),\omega_{f_{*}X} \big),
  \{ \Theta_{\tau} \} \big) \bigg] \in Z_n(D).$$
\end{construction}

Note that the class $[ (( (f_{*}X,|f_{*}X|,\{ \vartheta_\tau
\}),\omega_{f_{*}X}), \{ \Theta_{\tau} \}) ]$ is independent of
the choices we made. Thus construction
\ref{constr-pushfoftropcycles} immediately leads to the following

\begin{corollary}[Push-forward of tropical cycles] \label{coro-pushforwardoftropcycles}
  Let $C \in Z_n$ and $D \in Z_m$ be two cycles and let $f: C \rightarrow D$ be a morphism.
  Then for all $k$ there is a well-defined and $\ZZ$-linear map
  $$ Z_k(C) \longrightarrow Z_k(D): E \longmapsto f_{*}E :=
  (f|_E)_{*}E.$$
\end{corollary}
\begin{proof}
  The linearity can be proven similar to the affine case (cf.
  proposition \ref{prop-pushcycles}).
\end{proof}

Our next aim is to define the pull-back of Cartier divisors. But
first we need the following

\begin{lemma} \label{lemma-refinementofY}
  Let $C \in Z_n$ and $D \in Z_m$ be two tropical cycles and let
  $f: C \rightarrow D$ be a morphism. By definition there exist reduced
  representatives $\left(((X,|X|,\{\varphi_{\sigma}\}), \omega_X), \{\Phi_{\sigma}\} \right)$
  of $C$ and $\left(((Y,|Y|,\{\psi_{\tau}\}),\omega_Y),\{\Psi_{\tau}\} \right)$ of $D$ such that
  properties (a) and (b) in definition \ref{def-morphoftropcycle}
  are fulfilled. Let $\left(((Y_1,|Y_1|,\{\psi_{\tau'}'\}),\omega_{Y_1}),\{\Psi_{\tau'}'\} \right)$
  be a refinement of $Y$. Then there is a refinement
  $\left(((X_1,|X_1|,\{\varphi_{\sigma'}'\}), \omega_{X_1}), \{\Phi_{\sigma}'\} \right)$
  of $X$ such that properties (a) and (b) of definition \ref{def-morphoftropcycle}
  are fulfilled for $X_1$ and $Y_1$.
\end{lemma}
\begin{proof}
  Let $X_1 := \{ \sigma \cap f^{-1}(\tau)| \sigma \in X, \tau \in
  Y_1 \}$. By property (b) of definition \ref{def-morphoftropcycle}
  all $\varphi_{\sigma}(\sigma \cap f^{-1}(\tau))$ are rational
  polyhedra in $\RR^{n_{\sigma}}$. For every $\sigma' \in X_1$ choose $\sigma \in X$
  such that $\sigma' = \sigma \cap f^{-1}(\tau)$ for some $\tau \in Y_1$. Then we can define
  $\varphi_{\sigma'}' := \varphi_{\sigma}|_{\sigma'}:
  \sigma' \rightarrow H_{\sigma,\sigma'} \cong \RR^{n_{\sigma'}}$, where $H_{\sigma,\sigma'}$
  is the smallest affine subspace of $\RR^{n_{\sigma}}$ containing $\varphi_{\sigma}(\sigma')$.
  Moreover let $|X_1|:=|X|$. Note that with these settings
  $(X_1,|X_1|,\{ \varphi_{\sigma'}' \})$
  is a polyhedral complex. We can endow it with the weight function
  $\omega_{X_1}$ and the tropical structure
  $\{\Phi_{\sigma'}' \}$ induced by $X$. Together with $Y_1$ the tropical
  polyhedral complex $\left( ((X_1,|X_1|,\{ \varphi_{\sigma'}'\}),\omega_{X_1}),
  \{\Phi_{\sigma'}' \} \right)$ fulfills the requirements (a) and (b)
  of definition \ref{def-morphoftropcycle}.
\end{proof}

\begin{proposition}[Pull-back of Cartier divisors] \label{propo-pullbackofcartdivisors}
  Let $C \in Z_n$ and $D \in Z_m$ be tropical cycles and let $f: C
  \rightarrow D$ be a morphism. Then there is a well-defined and
  $\ZZ$-linear map
  $$\Div(D) \longrightarrow \Div(C): [\{(U_i,h_i)\}]
  \longmapsto f^{*}[\{(U_i,h_i)\}] := [\{(f^{-1}(U_i),h_i \circ f)\}].$$
\end{proposition}
\begin{proof}
  We have to show that $h \circ f \in \calK^{*}_C(f^{-1}(U))$
  for $h \in \calK^{*}_D(U)$ and that $h \circ f \in \calO^{*}_C(f^{-1}(U))$
  for $h \in \calO^{*}_D(U)$. Then the rest is obvious.\\
  So let $h \in \calK^{*}_D(U)$. Then there exists a representative
  $\left(((Y,|Y|,\{\psi_{\sigma}\}),\omega_Y),\{\Psi_{\tau}\} \right)$
  of $D$ such that for every polyhedron $\sigma \in Y$ the map $h \circ \psi_{\sigma}^{-1}$
  is locally integer affine linear. Moreover, since $f$ is a
  morphism there exist representatives $\left(((X,|X|,\{\varphi_{\sigma}\}),\omega_X),\{\Phi_{\tau}\} \right)$
  of $C$ and $\left(((Y',|Y'|,\{\psi_{\sigma'}'\}),\omega_{Y'}),\{\Psi_{\tau'}'\} \right)$
  of $D$ such that properties (a) and (b) of definition \ref{def-morphoftropcycle}
  are fulfilled, i.e. $f(\sigma) \subseteq \widetilde{\sigma} \in Y'$
  for all $\sigma \in X$ and the maps $\Psi_{\widetilde{\sigma}} \circ f \circ
  \Phi_{\sigma}^{-1}$ induce morphisms of fans. By lemma
  \ref{lemma-refinementofY} we may assume that $Y=Y'$. Now let
  $\sigma \in X$ and choose some $\widetilde{\sigma} \in Y$ such that
  $f(\sigma) \subseteq \widetilde{\sigma}$. Property (b) of
  definition \ref{def-morphoftropcycle} implies that
  $\psi_{\widetilde{\sigma}} \circ f \circ \varphi_{\sigma}^{-1}$
  and $\Psi_{\widetilde{\sigma}} \circ f \circ \Phi_{\sigma}^{-1}$
  are integer affine linear. Thus $h \circ f \circ
  \varphi_{\sigma}^{-1} = (h \circ \psi_{\widetilde{\sigma}}^{-1}) \circ (\psi_{\widetilde{\sigma}}
  \circ f \circ \varphi_{\sigma}^{-1})$ is locally integer affine
  linear and $h \circ f \in \calK^{*}_C(f^{-1}(U))$.
  If additionally $h \circ \Psi_{\widetilde{\sigma}}^{-1}$ is locally integer
  affine linear then so is $h \circ f \circ \Phi_{\sigma}^{-1} = (h
  \circ \Psi_{\widetilde{\sigma}}^{-1}) \circ
  (\Psi_{\widetilde{\sigma}} \circ f \circ \Phi_{\sigma}^{-1})$.
  Hence $h \circ f \in \calO^{*}_C(f^{-1}(U))$ for $h \in
  \calO^{*}_D(U)$.
\end{proof}

Our last step in this chapter is to state the analogon of the
projection formula from \ref{prop-projectionf}:

\begin{proposition}[Projection formula]\label{prop-abstrprojectionformula}
  Let $C \in Z_n$ and $D \in Z_m$ be two cycles and $f: C
  \rightarrow D $ be a morphism. Let $E \in Z_k(C)$ be a subcycle
  of $C$ and $d \in \Div(D)$ be a Cartier divisor. Then
  the following holds:
  $$d \cdot (f_{*}C) = f_{*}(f^{*}d \cdot C) \in Z_{k-1}(D).$$
\end{proposition}
\begin{proof}
  The claim follows from the constructions of $f_{*}C$ and $f^{*}d$, from definition
  \ref{def-intersectionproducts} and proposition \ref{prop-projectionf}.
\end{proof}

  \section{Rational Equivalence} \label{sec-rationalequivalence}

We will now make some first steps in establishing a concept of rational equivalence.

We fix an abstract tropical cycle $A$ as ambient space and an arbitrary subgroup $R \subseteq \Div(A)$ of the group of Cartier divisors on $A$. We define the \df{Picard group} as the quotient group $\Pic(A) := \Div(A) / R$. Let $R_k \subseteq Z_k(A)$ denote the group generated by $\{\varphi \cdot C | \varphi \in R, C \in Z_{k+1}(A)\}$, i.e. by all $k$-dimensional cycles obtained by intersecting a Cartier divisor from $R$ with an arbitrary $(k+1)$-dimensional cycle. We define the \df{$k$-th Chow group} to be $A_k(A) := Z_k(A) / R_k$.

\begin{corollary}[Intersection products modulo rational equivalence] \label{cor-intersectionsmodulorationalequivalence}
    The map
    \begin{eqnarray*}
        \cdot \; : \Pic(A) \times A_k(A) & \rightarrow & A_{k-1}(A), \\
        ([\varphi] , [D])      & \mapsto & [\varphi \cdot D]
    \end{eqnarray*}
    is well-defined and bilinear.
\end{corollary}

\begin{proof}
    By definition, for each $\varphi \in R$, $D \in Z_k(A)$ we have $\varphi \cdot D \in R_{k-1}$. Let furthermore $\varphi \cdot C$ be an element in $R_k$ (where $\varphi \in R, C \in Z_k(A)$). Then it follows from proposition \ref{prop-balancingconditionandcommutativity} b) that for arbitrary $\psi \in \Div(A)$ we get $\psi \cdot (\varphi \cdot C) = \varphi \cdot (\psi \cdot C) \in R_{k-1}$. The claim follows from the bilinearity of the intersection product.
\end{proof}

So far, our intersection theory takes place (at least locally) in $\RR^n$, which can be considered as the $n$-dimensional tropical algebraic torus. Especially, if we generated rational equivalence by all rational functions on $A$, the resulting Chow groups and intersection products would be useless in enumerative geometry: As in the classical case, the divisor of a rational function might have components in the ``boundary'' of some compactification of the ``affine'' variety $\RR^n$. Therefore, in the following we restrict the functions that generate rational equivalence to those ``whose divisor in any torical compactification has no components in the boundary''.

\begin{definition}[Rational equivalence generated by bounded functions] \label{def-boundedfunctions}
    Let $A$ be an abstract tropical cycle and $R(A) := \{[(|A|, \varphi)] | \varphi \text{ bounded} \}$ be the group of all Cartier divisors globally given by a bounded rational function. We define the \df{Picard group} $\Pic(A) := \Div(A) / R(A)$ and the \df{Chow groups} $A_k(A)$ as above. We call two Cartier divisors (two $k$-dimensional subcycles resp.) \df{rationally equivalent}, if their classes in $\Pic(A)$ ($A_k(A)$ resp.) are the same.
\end{definition}

Let us prove that we do not divide out too much for applications in enumerative geometry.

\begin{lemma} \label{lemma-degreezero}
    Let $C$ be an one-dimensional abstract tropical cycle, $\varphi \in R(C)$ a bounded rational function on $C$ and $(((X, |X|, \{m_\sigma\}_{\sigma \in X}), \omega_X), \{M_\sigma\}_{\sigma \in X})$ a representative of $C$ such that $|X| = |C|$ and for all $\sigma \in X$ the map $\varphi \circ m_\sigma^{-1} =: \varphi_\sigma$ is integer affine linear. Then
    $$
        \sum_{\{p\} \in X^{(0)}} \omega_\varphi (\{p\}) = 0,
    $$
    i.e. $\varphi \cdot C$ is of degree zero.
\end{lemma}

\begin{proof}
    By definition, for all $\{p\} \in X^{(0)}$ we have
    $$
        \omega_\varphi (\{p\}) = \sum_{\substack{\sigma \in X^{(1)} \\ p \in \sigma}} \omega(\sigma) \varphi_\sigma(u_{\sigma/\{p\}}).
    $$
    Note that if $\sigma \in X^{(1)}$ contains two different vertices, say $\partial_1 \sigma$ and $\partial_2 \sigma$, we have $u_{\sigma/\{\partial_1 \sigma\}} = - u_{\sigma/\{\partial_2 \sigma\}}$. If, otherwise, $\sigma$ contains less than two vertices, $m_\sigma(\sigma)$ is a non-compact polyhedron and therefore $\varphi$ can only be bounded if it is constant on $\sigma$. Together we get
    \begin{eqnarray*}
        \sum_{\{p\} \in X^{(0)}} \omega_\varphi (\{p\})
            & = &
            \sum_{\{p\} \in X^{(0)}} \sum_{\substack{\sigma \in X^{(1)} \\ p \in \sigma}}
                \omega(\sigma) \varphi_\sigma(u_{\sigma/\{p\}}) \\
        & = &
            \sum_{\substack{\sigma \in X^{(1)} \\ \exists ! \, \partial \sigma \, \in \sigma}}
                \omega(\sigma) \varphi_\sigma(u_{\sigma/\{\partial \sigma\}}) \\
        & &
            + \sum_{\substack{\sigma \in X^{(1)} \\ \exists ! \, \partial_1 \sigma,
                \partial_2 \sigma \, \in \sigma}}
                \omega(\sigma) \varphi_\sigma(u_{\sigma/\{\partial_1 \sigma\}})
                + \omega(\sigma) \varphi_\sigma(u_{\sigma/\{\partial_2 \sigma\}}) \\
        & = &
            \sum_{\substack{\sigma \in X^{(1)} \\ \exists ! \, \partial \sigma \, \in \sigma}}
                \omega(\sigma) \cdot 0 \\
        & &
            + \sum_{\substack{\sigma \in X^{(1)} \\ \exists ! \, \partial_1 \sigma,
                \partial_2 \sigma \, \in \sigma}}
                \omega(\sigma) \big(\underbrace{\varphi_\sigma(u_{\sigma/\{\partial_1 \sigma\}})
                - \varphi_\sigma(u_{\sigma/\{\partial_1 \sigma\}})}_{= 0}\big) \\
        & = & 0.
    \end{eqnarray*}
\end{proof}

\begin{remark} \label{redundantremark}
    As a consequence, for any cycle $C \in Z_*(A)$ there is a well-defined morphism
    $$
        \text{deg}: A_0(C) \longrightarrow \ZZ: [\lambda_1P_1+\ldots+\lambda_rP_r] \longmapsto
    \lambda_1+\ldots+\lambda_r.
    $$
    For $D \in A_0(C)$ the number $\text{deg}(D)$ is called the
    \df{degree} of $D$. \\
    Moreover, by corollary \ref{cor-intersectionsmodulorationalequivalence} there is a well-defined map of top products
    $$
        \Pic(A)^d \longrightarrow \ZZ:
            ([\varphi_1], \ldots, [\varphi_{\dim(C)}]) \longmapsto \text{deg}([\varphi_1 \cdot \ldots \cdot \varphi_{\dim(C)} \cdot C]),
    $$
        where $A$ is our ambient cycle and $d$ is the dimension
        of $C$. Of course, this map is of particular interest
        when dealing with enumerative questions.
\end{remark}

Of course, our chosen rational equivalence $R(A) := \{[(|A|, \varphi)] | \varphi \text{ bounded} \}$ should also be compatible with pull-back and push-forward. However, in the push-forward case we face problems due to our definition of rational functions. Let us first state the positive result in the pull-back case.

\begin{lemma}[Pull-back of rational equivalence] \label{lemma-pullbackratequiv}
    Let $C, D$ be tropical cycles and let $f : C \rightarrow D$ be a morphism between them. Then the pull-back map $\Div(D) \rightarrow \Div(C), \varphi \mapsto f^* \varphi$ induces a well-defined map on the quotients $\Pic(D) \rightarrow \Pic(C), [\varphi] \mapsto [f^* \varphi]$.
\end{lemma}
\begin{proof}
    We only have to show that for each element $(|D|, \psi) \in R(D)$ the pull-back Cartier divisor $f^* (|D|, \psi)$ lies in $R(C)$. But this follows from the trivial fact that the composition $\psi \circ f$ of a bounded function $\psi$ and an arbitrary map $f$ is again bounded.
\end{proof}

\begin{remark}[Push-forward of rational equivalence] \label{remark-pushforwardratequiv}
    The corresponding statement for push-forwards is false! Let us again consider the elliptic curve $C$ from Example \ref{example-cartierdivisor}. On this curve, the Weil divisor associated to the bounded rational function $\psi$ illustrated in the picture below equals $\divisor(\psi) = \alpha_1 + \alpha_2 - \alpha_3 - \alpha_4$.
    \medskip \\
    \begin{minipage}{\linewidth}
    \begin{center}
      \input{pic/RmPushFw.pstex_t}\\
    \end{center}
    \end{minipage}
    \medskip \\
    Let us now consider the cycle $D$ obtained by identifying $\alpha_1$ with $\alpha_3$ and the canonical projection map $f : C \rightarrow D$.
    \medskip \\
    \begin{minipage}{\linewidth}
    \begin{center}
      \input{pic/RmPushFw2.pstex_t}\\
    \end{center}
    \end{minipage}
    \medskip \\
    The push-forward of $\divisor(\psi)$ under this morphism is $f_* \divisor(\psi) = \alpha_2 - \alpha_4$. But this Weil divisor can obviously not be obtained by a rational function on $D$. This problem is due to our restrictive definition of rational functions (see remark \ref{remark-zerofunction}). We are currently working on a refined version of the related definitions.
\end{remark}

  \section{Intersection of cycles in $\RR^n$} \label{sec-stableintersection}

So far we are only able to intersect Cartier divisors with cycles.
Our aim in this section is now to define the intersection of two
cycles with ambient cycle $\RR^n$ (with trivial structure maps).
But first we need some preparations:

\begin{definition} \label{def-cartproductofcycles}
  Let $(((X,|X|,\{ \varphi_{\sigma} \}),\omega_X),\{\Phi_{\sigma}\})$
  and $(((Y,|Y|,\{ \psi_{\tau} \}),\omega_Y),\{\Psi_{\tau}\})$ be
  tropical polyhedral complexes. We denote by
  $$(((X,|X|,\{ \varphi_{\sigma} \}),\omega_X),\{\Phi_{\sigma}\})
  \times (((Y,|Y|,\{ \psi_{\tau} \}),\omega_Y),\{\Psi_{\tau}\})$$
  their \df{cartesian product}
  $$(((X \times Y,|X| \times |Y|,\{ \vartheta_{\sigma \times \tau}\}),\omega_{X \times Y}),\{\Theta_{\sigma \times \tau}\}),$$
  where
  \begin{eqnarray*}
    X \times Y & := & \left\{ \sigma \times \tau | \sigma \in X, \tau \in Y \right\},\\
    \vartheta_{\sigma \times \tau} & := & \varphi_{\sigma} \times \psi_{\tau} : \sigma \times \tau \longrightarrow \RR^{n_{\sigma}} \times \RR^{n_{\tau}},\\
    \omega_{X \times Y}(\sigma \times \tau) & := & \omega_X(\sigma) \cdot \omega_Y(\tau),\\
    \Theta_{\sigma \times \tau} & := & \Phi_{\sigma} \times \Psi_{\tau}: S_{\sigma}^X \times S_{\tau}^Y \longrightarrow |F_{\sigma}^X| \times |F_{\tau}^Y|.
  \end{eqnarray*}
  Let $\widetilde{F}_{\sigma}^X$ and $\widetilde{F}_{\tau}^Y$ be
  the entire fans associated with $F_{\sigma}^X$ and
  $F_{\tau}^Y$ from above. Obviously, the product ${\widetilde{F}_{\sigma}^X \times \widetilde{F}_{\tau}^Y}
  :=\{ \alpha \times \beta| \alpha \in \widetilde{F}_{\sigma}^X, \beta \in \widetilde{F}_{\tau}^Y \}$ with weight function
  ${\omega_{\widetilde{F}_{\sigma}^X \times \widetilde{F}_{\tau}^Y}(\alpha \times \beta)} := \omega_{\widetilde{F}_{\sigma}^X}(\alpha)
  \cdot \omega_{\widetilde{F}_{\tau}^Y}(\beta)$ is again a tropical fan and thus its intersection with $|F_{\sigma}^X| \times
  |F_{\tau}^Y|$ yields an open fan (cf. definition \ref{defn-openfan}). Hence the cartesian product
  ${(((X \times Y,|X| \times |Y|,\{ \vartheta_{\sigma \times \tau}\}),\omega_{X \times Y}),\{\Theta_{\sigma \times
  \tau}\})}$ is again a tropical polyhedral complex.\\

  If $C=[(X,\omega_X)]$ and $D=[(Y,\omega_Y)]$ are tropical cycles
  we define $$C \times D := [(X, \omega_X) \times (Y, \omega_Y)]$$
  for $(X, \omega_X) \times (Y, \omega_Y)$ as defined above. Note
  that $C \times D$ does not depend on the choice of the
  representatives $X$ and $Y$.
\end{definition}

\begin{remark} \label{rem-expressdiagonal}
We can express the diagonal in $\RR^n \times \RR^n$
$$[(\triangle,1)]=[(\{(x,x) | x \in \RR^n\},1)] \in Z_n(\RR^n
\times \RR^n)$$ as a product of Cartier divisors, namely
$$[(\triangle,1)] = \psi_1 \cdots \psi_n \cdot
\RR^n \times \RR^n,$$ where $\psi_i = [ \{
(\RR^n,\max\{0,x_i-y_i\}) \}] \in \Div(\RR^n \times \RR^n)$,
$i=1,\dots,n$. We will use this ability to define the intersection
product of any two cycles in $\RR^n$.
\end{remark}

\begin{definition} \label{def-intersectionproductRn}
  Let $\pi: \RR^n \times \RR^n \rightarrow \RR^n: (x,y) \mapsto
  x$. Then we define the intersection product of cycles in $\RR^n$
  by
  \begin{eqnarray*}
    Z_{n-k}(\RR^n) \times Z_{n-l}(\RR^n) & \longrightarrow & Z_{n-k-l}(\RR^n)\\
    (C,D) & \longmapsto & C \cdot D := \pi_{*}(\triangle \cdot (C
    \times D)),
  \end{eqnarray*}
  where $\pi_{*}$ denotes the push-forward as defined in \ref{coro-pushforwardoftropcycles}
  and $\triangle \cdot (C \times D) := \psi_1 \cdots \psi_n \cdot (C \times D)$  with
  $\psi_1,\dots,\psi_n$ as defined in remark \ref{rem-expressdiagonal}.
\end{definition}

Having defined this intersection product of arbitrary cycles in
$\RR^n$ we will prove now some basic properties. But as a start we
need the following lemmas:

\begin{lemma} \label{lemma-RntimesDequalsD}
  Let $C \in Z_k(\RR^n)$ be a cycle with representative $(X, \omega_X)$ and let $\psi_1,\dots,\psi_n$ be
  the Cartier divisors defined in remark \ref{rem-expressdiagonal}. Then
  $(X_j, \omega_{X_j})$ with
  $$X_j := \big\{ (\RR^n \times \sigma) \cap \{(x,y) \in \RR^n \times \RR^n| x_i=y_i \text{ for } i=j,\dots,n \} |\sigma \in X \big\},$$
  $$\omega_{X_j}\big((\RR^n \times \sigma) \cap \{(x,y) \in \RR^n \times \RR^n| x_i=y_i \text{ for } i=j,\dots,n \}\big) := \omega_X(\sigma)$$
  is a representative of $\psi_j \cdots \psi_n \cdot \RR^n \times C$.
\end{lemma}
\begin{proof}
  We use induction on $j$. For $j=n+1$ there is nothing to show.
  Now let the above representative be correct for some $j+1$. We
  have to show that $X_{j}$ is a tropical polyhedral complex and
  that it represents $\psi_{j} \cdots \psi_n \cdot \RR^n \times
  C$: Note that
  \[ \begin{array}{c}
    \dim \left( (\RR^n \times \sigma) \cap \{(x,y) \in \RR^n \times \RR^n| x_i=y_i \text{ for } i=j,\dots,n \}\right) \\
    < \dim \left( (\RR^n \times \sigma) \cap \{(x,y) \in \RR^n \times \RR^n| x_i=y_i \text{ for } i=j+1,\dots,n \}\right)
  \tag{$*$} \label{tag-diminc} \end{array} \]
  for all $\sigma \in X$. Hence $X_{j}$ is a tropical polyhedral
  complex. Moreover note that
  $$ \widetilde{X}_{j+1} := \{ \sigma \cap \{x_{j}-y_{j}=0\}, \sigma \cap \{x_{j}-y_{j} \leq
  0\}, \sigma \cap \{x_{j}-y_{j} \geq 0\} | \sigma \in X_{j+1}\}$$
  with weights induced by $X_{j+1}$ is a refinement of $X_{j+1}$
  such that $\max\{0,x_j-y_j\}$ is linear on every face of
  $\widetilde{X}_{j+1}$. By (\ref{tag-diminc}) there are exactly
  two types of faces of codimension one in $\widetilde{X}_{j+1}$:
  \begin{enumerate}
    \item[(i)] $(\RR^n \times \sigma) \cap \{x_i-y_i=0 \text{ for } i=j,\dots,n
    \}$ with $\sigma \in X$, $\codim(\sigma) = 0$,
    \item[(ii)] $(\RR^n \times \sigma) \cap \{x_i-y_i=0 \text{ for }
    i=j+1,\dots,n; \enspace x_j-y_j \leq 0 \}$ or\\
    $(\RR^n \times \sigma) \cap \{x_i-y_i=0 \text{ for }
    i=j+1,\dots,n; \enspace x_j-y_j \geq 0 \}$with $\sigma \in X$,\\ $\codim(\sigma) = 1$,
  \end{enumerate}
  where the faces of the second type are not contained in
  $\{(x,y) \in \RR^n \times \RR^n| x_j = y_j \}$. Hence
  $\max\{0,x_j-y_j\}$ is linear on a neighborhood of every face of
  type (ii) and thus these faces get weight zero in
  $\max\{0,x_j-y_j\} \cdot \widetilde{X}_{j+1}$. The faces of type
  (i) are weighted by $\omega_{X_{j+1}}((\RR^n \times \sigma) \cap \{x_i-y_i=0 \text{ for } i=j+1,\dots,n
  \})$ in $\max\{0,x_j-y_j\} \cdot \widetilde{X}_{j+1}$ since
  ${x_1-y_1},\dots,x_n-y_n$ are part of a lattice basis of $(\ZZ^n
  \times \ZZ^n)^{\vee}$. Thus ${\max\{0,x_j-y_j\}} \cdot
  \widetilde{X}_{j+1}= X_j$ and $X_j$ is a representative of
  $\psi_j \cdots \psi_n \cdot \RR^n \times C$.
\end{proof}

\begin{corollary} \label{coro-RntimesCequalsC}
  Let $C \in Z_k(\RR^n)$ be a cycle. Then we have the equation:
  $$\RR^n \cdot C = C.$$
\end{corollary}
\begin{proof}
  Let $(X, \omega_X)$ be a representative of $C$, let $\pi: \RR^n \times \RR^n \rightarrow \RR^n:
  (x,y) \mapsto x$ and let $\psi_1, \dots, \psi_n$ be
  the Cartier divisors defined in remark \ref{rem-expressdiagonal}. By lemma
  \ref{lemma-RntimesDequalsD} we know that $X_1=\{ \{(x,x)|x \in \sigma\} | \sigma \in X \}$
  with $\omega_{X_1}(\{(x,x)|x \in \sigma\})=\omega_X(\sigma)$ is a representative of
  $\psi_1 \cdots \psi_n \cdot \RR^n \times C$. Hence $$\RR^n \cdot C =
  \pi_{*}(\psi_1 \cdots \psi_n \cdot \RR^n \times C) = [\pi_{*}(X_1,\omega_{X_1})] = [(X,\omega_X)] = C.$$
\end{proof}

\begin{lemma} \label{lemma-pullbackofcartproduct}
  Let $C \in Z_k(\RR^n)$ and $D \in Z_l(\RR^m)$ be abstract cycles, $\varphi \in \Div(\RR^n)$ a
  Cartier divisor and $\pi: {\RR^n \times \RR^m \rightarrow
  \RR^n:} (x,y) \mapsto x$. Then:
  $$(\varphi \cdot C) \times D= \pi^{*}\varphi \cdot (C \times D).$$
\end{lemma}
\begin{proof}
  We prove the statement for affine cycles $C,D$ and an affine
  Cartier divisor $\varphi$. The general case then follows by
  applying the statement locally.\\
  Choose arbitrary representatives $Y$ of $D$ and $h$ of
  $\varphi$ and choose a representative $X$ of $C$ such that
  $h$ is linear on every face of $X$. This implies that $\pi^{*}h$
  is linear on every face of $X \times Y$, too. In $X \times Y$ we have
  two types of faces of codimension one:
  \begin{enumerate}
    \item[(i)] $\sigma \times \tau$ with $\sigma \in X, \tau \in
    Y, \codim(\sigma)=1, \codim(\tau)=0$,
    \item[(ii)] $\sigma \times \tau$ with $\sigma \in X, \tau \in
    Y, \codim(\sigma)=0, \codim(\tau)=1$.
  \end{enumerate}
  For the second type the adjacent facets are exactly all $\sigma
  \times \widetilde{\tau}$ with $\widetilde{\tau} > \tau$. We get
  $\omega_h(\sigma \times \tau)=0$ in $h \cdot X \times Y$ as $\pi^{*}h$ is linear on
  $\sigma \times |Y|$. For the first type the adjacent facets are exactly all $\widetilde{\sigma}
  \times \tau$ with $\widetilde{\sigma} > \sigma$ and the weights
  can be calculated exactly like for $h \cdot X$. This finishes
  the proof.
\end{proof}

Let $C$ and $D$ be cycles in $\RR^n$. Assume that $C$ can be
expressed as a product of Cartier divisors, i.e. there are
$\varphi_1,\dots,\varphi_r \in \Div(\RR^n)$ such that $C =
\varphi_r \cdots \varphi_1 \cdot \RR^n$. The obvious questions are
now how $C \cdot D$ relates to $\varphi_r \cdots \varphi_1 \cdot
D$ and whether $\varphi_r \cdots \varphi_1 \cdot D$ depends on the
choice of the Cartier divisors $\varphi_i$. To answer this
question we first prove a somewhat stronger statement:

\begin{lemma} \label{lemma-cartdivinterchanges}
  Let $C \in Z_k(\RR^n)$ and $D \in Z_l(\RR^n)$ be cycles and $\varphi \in
  \Div(\RR^n)$ a Cartier divisor. Then we have the equality:
  $$(\varphi \cdot C) \cdot D = \varphi \cdot (C \cdot D).$$
\end{lemma}
\begin{proof}
  Let $\pi: \RR^n \times \RR^n \rightarrow \RR^n: (x,y) \mapsto
  x$ be like above. It holds:
  \begin{eqnarray*}
    (\varphi \cdot C) \cdot D &=& \pi_{*}(\triangle \cdot (\varphi \cdot C) \times D)\\
    &\stackrel{\ref{lemma-pullbackofcartproduct}}{=}& \pi_{*}(\pi^{*}\varphi \cdot \triangle \cdot C \times D)\\
    &\stackrel{\ref{prop-abstrprojectionformula}}{=}& \varphi \cdot \pi_{*}(\triangle \cdot C \times D)\\
    &=& \varphi \cdot (C \cdot D).
  \end{eqnarray*}
\end{proof}

\begin{corollary} \label{coro-relationprodwithcartdiv}
  Let $C \in Z_k(\RR^n)$ be a cycle such that there are Cartier divisors \linebreak
  $\varphi_1,\dots,\varphi_r \in \Div(\RR^n)$ with $\varphi_r \cdots \varphi_1 \cdot \RR^n = C$
  and let $D \in Z_l(\RR^n)$ be any cycle.
  Then $$\varphi_r \cdots \varphi_1 \cdot D = C \cdot D.$$
\end{corollary}
\begin{proof}
  Applying lemma \ref{lemma-cartdivinterchanges} and lemma \ref{lemma-RntimesDequalsD} we obtain
  $$ C \cdot D = (\varphi_r \cdots \varphi_1 \cdot \RR^n) \cdot D = \varphi_r \cdots \varphi_1 \cdot (\RR^n \cdot
  D) = \varphi_r \cdots \varphi_1 \cdot D.$$
\end{proof}

\begin{remark} \label{rem-indepenceofchoiceofdivisors}
  Note that corollary \ref{coro-relationprodwithcartdiv} in particular
  implies that our definition of the intersection product on $\RR^n$
  (cf. \ref{def-intersectionproductRn}) is independent of the choice of the
  Cartier divisors describing the diagonal $\triangle$.
\end{remark}

\begin{theorem} \label{thm-propertiesofintproduct}
  Let $C, C' \in Z_k(\RR^n)$, $D \in Z_l(\RR^n)$ and $E \in Z_m(\RR^n)$ be cycles.
  Then the following equations hold:
  \begin{enumerate}
    \item $C \cdot D= D \cdot C$,
    \item $(C + C') \cdot D = C \cdot D + C' \cdot D$,
    \item $(C \cdot D) \cdot E = C \cdot (D \cdot E)$.
  \end{enumerate}
\end{theorem}
\begin{proof}
  (a): Let $\psi_1, \dots, \psi_n \in \Div(\RR^n \times \RR^n)$ be
  like defined in remark \ref{rem-expressdiagonal}. Note that for
  every $i \in \{1,\dots,n\}$ the maps $\max\{0,x_i-y_i\}$ and $\max\{0,y_i-x_i\}$ only
  differ by a globally linear map and hence define the same Cartier divisor.
  Thus we get $$\pi_{*}(\psi_1 \cdots \psi_n
  \cdot C \times D) = \pi_{*}(\psi_1 \cdots \psi_n \cdot D \times C).$$
  (b): Follows immediately by bilinearity of the intersection
  product $$\Div(\RR^n \times \RR^n) \times Z_{p}(\RR^n \times
  \RR^n) \stackrel{\cdot}{\longrightarrow} Z_{p-1}(\RR^n \times
  \RR^n),$$
  linearity of the push-forward and the fact that $(C + C') \times D = C \times D + C' \times D$.\\
  (c): We will show that $\triangle \cdot C \times (\pi_{*}(\triangle \cdot D \times E)) =
  \triangle \cdot (\pi_{*}(\triangle \cdot C \times D) \times E):$\\
  Let $\pi^{12}: (\RR^n)^3 \rightarrow (\RR^n)^2: (x,y,z) \mapsto (x,y)$,
  $\pi^{13}: (\RR^n)^3 \rightarrow (\RR^n)^2: (x,y,z) \mapsto (x,z)$
  and $\pi^{23}: (\RR^n)^3 \rightarrow (\RR^n)^2: (x,y,z) \mapsto (y,z)$.
  An easy calculation shows that
  \begin{eqnarray} \triangle \cdot C \times (\pi_{*}(\triangle \cdot D \times E))
  = \triangle \cdot \pi_{*}^{12}(C \times (\triangle \cdot D \times
  E))\end{eqnarray}
  and
  \begin{eqnarray} \triangle \cdot (\pi_{*}(\triangle \cdot C \times D) \times E)
  = \triangle \cdot \pi_{*}^{13}((\triangle \cdot C \times D) \times
  E).\end{eqnarray}
  Now let $\psi_1, \dots, \psi_n $ be the Cartier divisors
  defined in remark \ref{rem-expressdiagonal}. We label these
  Cartier divisors with pairs of letters $\psi_i^{xy}$ to point out
  the coordinates they are acting on. We obtain
  \renewcommand{\arraystretch}{1.6}
  $$\begin{array}{l}
    \triangle \cdot C \times (\pi_{*}(\triangle \cdot D \times E))\\
    \stackrel{(1)}{=} \enspace \triangle \cdot \pi_{*}^{12}(C \times (\triangle \cdot D \times E)) \\
    = \enspace \psi_1^{xy} \cdots \psi_n^{xy} \cdot \pi_{*}^{12}(C \times (\psi_1^{yz} \cdots \psi_n^{yz} \cdot D \times E)) \\
    \stackrel{\ref{prop-abstrprojectionformula}}{=} \enspace
    \pi_{*}^{12}((\pi^{12})^{*}\psi_1^{xy} \cdots (\pi^{12})^{*}\psi_n^{xy} \cdot C \times (\psi_1^{yz} \cdots
    \psi_n^{yz} \cdot D \times E))\\
    \stackrel{\ref{lemma-pullbackofcartproduct}}{=} \enspace
    \pi_{*}^{12}((\pi^{23})^{*}\psi_1^{yz} \cdots (\pi^{23})^{*}\psi_n^{yz} \cdot (\pi^{12})^{*}\psi_1^{xy} \cdots (\pi^{12})^{*}\psi_n^{xy} \cdot C \times D \times E)\\
    \stackrel{\ref{coro-relationprodwithcartdiv}}{=} \enspace
    \pi_{*}^{13}((\pi^{12})^{*}\psi_1^{xy} \cdots (\pi^{12})^{*}\psi_n^{xy} \cdot (\pi^{13})^{*}\psi_1^{xz} \cdots (\pi^{13})^{*}\psi_n^{xz} \cdot C \times D \times E)\\
    \stackrel{\ref{lemma-pullbackofcartproduct}}{=} \enspace
    \pi_{*}^{13}((\pi^{13})^{*}\psi_1^{xz} \cdots (\pi^{13})^{*}\psi_n^{xz} \cdot (\psi_1^{xy} \cdots \psi_n^{xy} \cdot C \times D) \times E)\\
    \stackrel{\ref{prop-abstrprojectionformula}}{=} \enspace
    \psi_1^{xz} \cdots \psi_n^{xz} \cdot \pi_{*}^{13}((\psi_1^{xy} \cdots \psi_n^{xy} \cdot C \times D) \times E)\\
    = \enspace
    \triangle \cdot \pi_{*}^{13}((\triangle \cdot C \times D) \times E)\\
    \stackrel{(2)}{=} \enspace
    \triangle \cdot (\pi_{*}(\triangle \cdot C \times D) \times E).
  \end{array}$$
  \renewcommand{\arraystretch}{1}
  This proves (d).
\end{proof}

It remains to show that our intersection product is well-defined
modulo rational equivalence. If this is the case the intersection
product induced on $A_{*}(\RR^n)$ clearly inherits the properties
of the intersection product on $Z_{*}(\RR^n)$ we have proven in
this section.

\begin{proposition} \label{prop-intprodRnratequiv}
  The intersection product $Z_{n-k}(\RR^n) \times Z_{n-l}(\RR^n)
  \stackrel{\cdot}{\longrightarrow} Z_{n-k-l}(\RR^n)$ induces a
  well-defined and bilinear map
  $$ A_{n-k}(\RR^n) \times A_{n-l}(\RR^n)
  \stackrel{\cdot}{\longrightarrow} A_{n-k-l}(\RR^n): ([C],[D])
  \longmapsto [C] \cdot [D] := [C \cdot D].$$
\end{proposition}
\begin{proof}
  Let $h \cdot C \in R_{n-k}$ (cf. section \ref{sec-rationalequivalence}) and $D \in Z_{n-l}(\RR^n)$.
  Using lemma \ref{lemma-cartdivinterchanges} we
  can conclude that $(h \cdot C) \cdot D = h \cdot (C \cdot D) \in R_{n-k-l}$.
\end{proof}

Our last step in this section is to prove a Bézout-style theorem
for a special class of tropical cycles in $\RR^n$ called
\df{$\PP^n$-generic cycles}. But first we need some further
definitions:

\begin{definition} \label{def-translationalongv}
  Let $X$ be a tropical polyhedral complex in $\RR^n$ and let $v \in \RR^n$.
  We denote by $X(v)$ the translation $$ X(v) := \{\sigma + v| \sigma \in X \}$$
  of $X$ along $v$. If $[X]=C \in Z_k(\RR^n)$ then $C(v):=[X(v)]$.
  Note that the class $C(v)$ is independent of the representative
  $X$.
\end{definition}

\begin{definition} \label{def-degreeofcycle}
  Let $C \in Z_k(\RR^n)$ be a tropical cycle and let $L^n_k$ be
  the tropical fan defined in example \ref{example-selfintersectionofhyperplanes}.
  Then we define the \df{degree of $C$} to be the number
  $$\deg(C) := \deg(C \cdot [L^n_{\codim X}]),$$ where the second
  map $\deg : Z_0(\RR^n) \rightarrow \ZZ: \lambda_1 P_1 + \ldots + \lambda_r P_r \mapsto \lambda_1+ \ldots+ \lambda_r$
  is the usual degree map. Then the map $\deg: Z_k(\RR^n) \rightarrow \ZZ$ is obviously linear by definition.
  Moreover, we define the degree of $[C] \in A_k(\RR^n)$
  to be $\deg([C]) := \deg(C).$ Note that $\deg([C])$ is
  well-defined by remark \ref{redundantremark}.
\end{definition}

\begin{lemma} \label{lemma-invarianceofdegree}
  Let $C \in Z_k(\RR^n)$ and $D \in Z_{n-k}(\RR^n)$ be two
  tropical cycles of complementary dimensions. Then $$\deg(C \cdot D)=\deg(C(v_1) \cdot
  D(v_2))$$ for all vectors $v_1,v_2 \in \RR^n$. In particular
  $\deg(C) = \deg(C(v))$ for all $v \in \RR^n$.
\end{lemma}
\begin{proof}
  Let $\pi:\RR^n \times \RR^n \rightarrow \RR^n: (x,y) \mapsto x$
  be the projection map as above and for $u=(u_1,\ldots,u_n) \in \RR^n$ let $$\triangle(u) \cdot (C \times D) :=
  \psi_1(u_1) \cdots \psi_n(u_n) \cdot (C \times D)$$ with
  $\psi_i(u_i) := [ \{(\RR^n,\max\{0,x_i-y_i+u_i\}) \}] \in \Div(\RR^n \times
  \RR^n)$ be the intersection with the translated diagonal (cf. definition \ref{def-intersectionproductRn}).
  Note that the rational function ${\max\{0,x_i-y_i\}} - {\max\{0,x_i-y_i+u_i\}}$ is bounded
  and that hence $[\psi_i]=[\psi_i(u_i)] \in \Pic(\RR^n \times \RR^n)$ for all ~$i$. It follows that
  $$[\triangle \cdot (C \times D)] = [\triangle(u) \cdot (C \times D)] \in A_0(\RR^n \times \RR^n)$$
  and thus we get
  \begin{eqnarray*}
    \deg(C \cdot D) &=& \deg(\pi_{*}( \triangle \cdot (C \times D)))\\
    &=& \deg(\triangle \cdot (C \times D))\\
    &\stackrel{\ref{redundantremark}}{=}& \deg(\triangle(v_1-v_2) \cdot (C \times D))\\
    &=& \deg(\triangle \cdot (C(v_1) \times D(v_2)))\\
    &=& \deg(\pi_{*}( \triangle \cdot (C(v_1) \times D(v_2))))\\
    &=& \deg(C(v_1) \cdot D(v_2)).
  \end{eqnarray*}
\end{proof}

\begin{definition}[$\PP^n$-generic cycles] \label{def-p2generic}
  Let $C \in Z_k(\RR^n)$ be a tropical cycle. $C$ is called
  \df{$\PP^n$-generic} if for one (and thus for every) representative
  $X$ of $C$ holds:
  For every face $\sigma \in X^{(k)}$ there exists a polytope
  $P_\sigma \subseteq \RR^n$ of some dimension $r \in \{0,\dots,k\}$
  and a cone $\widetilde{\sigma} \in (L^n_k)^{(k-r)}$
  such that $\sigma \subseteq P_\sigma + \widetilde{\sigma}$.
\end{definition}

\begin{theorem}[Bézout's theorem] \label{thm-bezout}
  Let $C \in Z_k(\RR^n)$ and $D \in Z_{n-k}(\RR^n)$ be two
  tropical cycles of complementary dimensions. Moreover, assume
  that $C$ and $D$ are $\PP^n$-generic. Then:
  $$ \deg(C \cdot D) = \deg(C) \cdot \deg(D).$$
\end{theorem}
\begin{proof}
  Let $(X,\omega_X)$ be a representative of $C$ and $(Y,\omega_Y)$
  be a representative of $D$. Moving $X$ along a (generic) direction vector $a=(a_1,\ldots,a_n) \in
  \RR^{k}_{\gg 0} \times \RR^{n-k}_{\ll 0}$ we can reach that $|X(a)|$ and
  $|Y|$ intersect in points in the interior of maximal faces only,
  namely $|X(a)| \cap |Y| = \{P_{ij}| i=1,\ldots,r; j=1,\ldots,s\}$
  with $P_{ij}= \sigma_i \cap \sigma'_j$ for facets (we use the
  notation introduced in example \ref{example-selfintersectionofhyperplanes}
  for the cones of $L_k^n$ here)
  \begin{itemize}
    \item $\sigma_i \in X(a)^{(k)} \text{ with } \sigma_i \subseteq
    \sigma_{\{1,\ldots,k\}} + u_i \in L^n_k(u_i) \text{ and}$
    \item $\sigma'_j \in Y^{(n-k)} \text{ with } \sigma'_j \subseteq
    \sigma_{\{k+1,\ldots,n\}} + v_j \in L^n_{n-k}(v_j)$.
  \end{itemize}
  Hence we can conclude that $X(a) \cdot Y = \sum_{i=1}^r \sum_{j=1}^s
  \omega_X(\sigma_i) \omega_Y(\sigma'_j) P_{ij}$ and thus by lemma \ref{lemma-invarianceofdegree} $${\deg(X
  \cdot Y)} = \deg(X(a) \cdot Y) = \sum_{i=1}^r \sum_{j=1}^s \omega_X(\sigma_i)
  \omega_Y(\sigma'_j).$$ Moreover we can deduce that $|X(a)| \cap |L^n_{n-k}(v_1)| =
  \{P_{11},\ldots,P_{r1}\}$. Hence $X(a) \cdot L^n_{n-k}(v_1) =
  \sum_{i=1}^r \omega_X(\sigma_i) P_{i1}$ and again by lemma \ref{lemma-invarianceofdegree} $$\deg(X) =
  \deg(X(a) \cdot L^n_{n-k}(v_1)) = \sum_{i=1}^r
  \omega_X(\sigma_i).$$
  Analogously we obtain $$\deg(Y) =
  \deg(Y \cdot L^n_{k}(u_1)) = \sum_{j=1}^s
  \omega_Y(\sigma'_j).$$
  Thus the claim follows.
  \begin{figure}
    \input{pic/ExBezout.pstex_t}\\
    The intersection of $X(a)$ and $Y$ as described in
    \ref{thm-bezout}.
  \end{figure}
\end{proof}

\end{document}